\documentclass[12pt]{article}
    \usepackage[margin=1in]{geometry}
    \usepackage[utf8]{inputenc}
    \usepackage[T1]{fontenc}    
    \usepackage{amsmath,amssymb,amsfonts,amsthm}
    \usepackage{mathtools}
    \usepackage{physics}
    \usepackage[symbol]{footmisc}
    \usepackage{url}
    \usepackage{bigints}
    \usepackage{graphicx}
    \usepackage{tabularx}
    \usepackage{xcolor}
    \usepackage{todonotes}
    \usepackage{etoc}
    \usepackage{chngcntr}    
\newtheorem{definition}{Defn.}
\newtheorem{cor}{Cor.}
\newtheorem{theo}{Theo.}

\graphicspath{{./Fig/}}

\begin{document}

\begin{titlepage}
\begin{center}

    \Large\textbf{Framing local structural identifiability and observability in terms of parameter-state symmetries}
       
    \normalsize
        
    \vspace{0.25cm}
    \setcounter{footnote}{0}
    \setlength{\footnotemargin}{0.8em}
    {\normalsize Johannes G Borgqvist\footnote{Corresponding author: \url{johborgq@chalmers.se}}\footnote{Mathematical Sciences, Chalmers University of Technology, Gothenburg, Sweden}, Alexander P Browning\footnote{Mathematical Institute, University of Oxford, United Kingdom\label{Oxford}}\footnote{School of Mathematics and Statistics, University of Melbourne, Melbourne, Australia\label{Melbourne}}, Fredrik Ohlsson\footnote{Department of Mathematics and Mathematical Statistics, Umeå University, Umeå, Sweden}, Ruth E Baker\footref{Oxford}}
    \setlength{\footnotemargin}{1.8em}
        
\abstract{We introduce a subclass of Lie symmetries, called parameter–state symmetries, to analyse the local structural identifiability and observability of mechanistic models consisting of state-dependent ODEs with observed outputs. These symmetries act on parameters and states while preserving observed outputs at every time point. We prove that locally structurally identifiable parameter combinations and locally structurally observable states correspond to universal invariants of all parameter–state symmetries of a given model. We illustrate the framework on four previously studied mechanistic models, confirming known identifiability results and revealing novel insights into which states are observable, providing a unified symmetry-based approach for analysing structural properties of dynamical systems.}    
    
\end{center}
\end{titlepage}

\section{Introduction}

Structural identifiability and structural observability are fundamental concepts in the analysis of mechanistic models, particularly those described by ordinary differential equations (ODEs). Structural identifiability asks whether the parameters of a model can be uniquely determined from perfect knowledge of the outputs, while structural observability asks whether the states of the model can be uniquely reconstructed from the outputs. Understanding these properties is crucial for meaningful parameter estimation and state reconstruction in mechanistic modelling.

The standard differential algebra approach re-writes the original state-dependent ODE system together with output equations as an output system by eliminating the states, producing a potentially higher-order system of ODEs that depend only on the observed outputs, the independent variable, and the parameters. A key mathematical property of this reformulation is that the original ODE system with outputs is \textit{output-equivalent} with the derived output system: if the original system is solved and the observed outputs generated, these same observed outputs also solve the output system~\cite{eisenberg2019inputoutputequivalenceidentifiabilitysimple}. This equivalence guarantees that the conclusions of the standard differential algebra approach provide valid information about the structural identifiability of the original state-dependent system~\cite{ljung1994global,walter1982global,hong2020global,eisenberg2019inputoutputequivalenceidentifiabilitysimple,renardy2022structural,dong2023structuralIdentifiabilityJL}. The standard differential algebra approach is primarily used to analyse global structural identifiability, whereas alternative methods such as the EAR method~\cite{dong2023structuralIdentifiabilityJL} provide a complementary framework for analysing local structural identifiability and structural observability.

An alternative approach to study local identifiability and observability is based on an extended form of Lie symmetries, referred to as \textit{full symmetries}, which act on both the states and parameters of mechanistic models. Previous works~\cite{massonis2020finding,castro2020structuralIdentifiability,yates2009structural,merkt2015higher,villaverde2022symmetries} explored the use of full symmetries to analyse local identifiability and observability, but the exact link between these full symmetries, on the one hand, and locally structurally identifiable parameter combinations and locally structurally observable state combinations, on the other, remained unclear. In our previous work~\cite{borgqvist2024framingglobalstructuralidentifiability}, we partially clarified this link by showing that local structural identifiability can be understood in terms of universal parameter invariants associated with so-called \textit{parameter symmetries} of the output system.

Building on this framework, we demonstrate that the universal parameter invariants in our previous work are exactly the same as the locally structurally identifiable parameter combinations identified by the \textit{parameter-state symmetries} introduced in this work. Better still, the extension introduced here which is a subclass of full symmetries—parameter-state symmetries—simultaneously transforms states and parameters in such a way that the \textit{observed outputs are preserved at each time point}, revealing the locally structurally observable state combinations in addition to the locally structurally identifiable parameter combinations. This formulation is novel: the structural identifiability part recovers the results of our previous work on parameter symmetries of the output system, while the structural observability part provides a new methodology for identifying locally observable state combinations.

To illustrate how these local structural properties can be analysed using parameter-state symmetries, we analyse a glucose-insulin regulation model~\cite{cobelli1980parameter,massonis2020finding} and an epidemiological SEI model of tuberculosis outbreaks~\cite{renardy2022structural}. These examples demonstrate how parameter-state symmetries can be used to study both local structural identifiability and observability simultaneously. For the interested reader, we provide detailed information about structural identifiability and observability as well as Lie symmetries in the Mathematical Preliminaries in Section~\ref{sec:math_preliminaries}. Readers primarily interested in applications are advised to move directly to Section~\ref{sec:results}.

In our previous work~\cite{borgqvist2024framingglobalstructuralidentifiability}, local structural identifiability was analysed using parameter symmetries of the output system. This work, which both builds upon and generalises those previous findings, provides a method for analysing local structural identifiability and observability simultaneously by analysing the full model consisting of a state-dependent ODE system and equations for observed outputs. Together, our parameter-state symmetries constitute a comprehensive framework for the systematic analysis of local structural identifiability and observability of mechanistic models.

\section{Mathematical Preliminaries}
\label{sec:math_preliminaries}

In this section, we summarise three key concepts that form the foundation of our symmetry-based approach for analysing local parameter identifiability and state observability. Our goal here is to provide a concise and accessible introduction to these tools, giving the reader a brief, self-contained overview of the concepts needed for our approach. First, we introduce definitions and notation for structural identifiability and structural observability, describing how parameters can be inferred from outputs and how states can be reconstructed from outputs~\cite{renardy2022structural,meshkat2014identifiability}. Second, we review a powerful tool in the analysis of ODEs, namely classical Lie symmetries, which are smooth transformations of independent (e.g. time $t$) and dependent variables (e.g. states $\mathbf{x}$) that map solutions to other solutions~\cite{bluman1989symmetries,olver2000applications,hydon2000symmetry}. Third, we discuss extended symmetries that also act on parameters in addition to states, that we refer to as full symmetries~\cite{massonis2020finding,villaverde2022symmetries,merkt2015higher}, which provide a useful conceptual bridge to the parameter-state symmetries introduced in this work.

\subsection{Mechanistic models consisting of ODE systems with outputs}
\label{sec:state_dep_ODE}
We consider a general mechanistic model consisting of a system of ordinary differential equations (ODEs) describing $n$ states $\mathbf{x} = (x_1,\ldots,x_n)^\top \in \mathbb{R}^n$ depending on time $t\in\mathbb{R}$ and $p$ parameters $\boldsymbol{\theta} = (\theta_1,\ldots,\theta_p)^\top \in \mathbb{R}^p$, together with $m$ observed outputs $\mathbf{y} = (y_1,\ldots,y_m)^\top$:

\begin{equation}
    \dfrac{\mathrm{d}\mathbf{x}}{\mathrm{d}t} = \mathbf{f}(t,\mathbf{x},\boldsymbol{\theta}), 
    \qquad 
    \mathbf{y} = \mathbf{h}(t,\mathbf{x},\boldsymbol{\theta}),
    \label{eq:ODE_sys_full}
\end{equation}
where $\mathbf{f}:\mathbb{R}\times\mathbb{R}^n\times\mathbb{R}^p \rightarrow \mathbb{R}^n$ and $\mathbf{h}:\mathbb{R}\times\mathbb{R}^n\times\mathbb{R}^p \rightarrow \mathbb{R}^m$ are sufficiently smooth functions. 

The model in Eq.~\eqref{eq:ODE_sys_full} consisting of a state-dependent ODE system together with equations for observed ouputs and its corresponding \textit{output system}, in which the states are eliminated and only the outputs $\mathbf{y}$, their time derivatives, time $t$, and the parameters $\boldsymbol{\theta}$ remain, are \emph{output-equivalent}~\cite{renardy2022structural,ljung1994global}. That is, any output trajectory generated by Eq.~\eqref{eq:ODE_sys_full} also solves the output system. Previous analyses of structural identifiability using the standard differential algebra approach focused on the output system~\cite{renardy2022structural,ljung1994global}, which is also the system analysed in our previous work on parameter symmetries~\cite{borgqvist2024framingglobalstructuralidentifiability}. In this work, by contrast, we analyse two structural properties of the full model in Eq.~\eqref{eq:ODE_sys_full}, namely \textit{local structural identifiability} and \textit{local structural observability}. 

To prepare for the local definitions of structural identifiability and observability, we introduce neighbourhoods around given states and parameter values. For a centre point $\mathbf{x}_C \in \mathbb{R}^n$ and radius $\rho_{\mathbf{x}}>0$, the ball in state space is defined as
\begin{equation}
\mathcal{B}_{\mathbf{x}}(\mathbf{x}_C) =
\left\{\mathbf{x} \in \mathbb{R}^n :
\|\mathbf{x}-\mathbf{x}_C\|_2 < \rho_{\mathbf{x}}\right\}.
\end{equation}
Similarly, for parameters with centre $\boldsymbol{\theta}_C \in \mathbb{R}^p$ and radius $\rho_{\boldsymbol{\theta}}>0$, the ball in parameter space is defined as
\begin{equation}
\mathcal{B}_{\boldsymbol{\theta}}(\boldsymbol{\theta}_C) =
\left\{\boldsymbol{\theta} \in \mathbb{R}^p :
\|\boldsymbol{\theta}-\boldsymbol{\theta}_C\|_2 < \rho_{\boldsymbol{\theta}}\right\}.
\end{equation}
These balls define neighbourhoods in state and parameter space within which local properties of the model can be analysed. In particular, they provide the natural setting for the local definitions of structural identifiability and observability introduced in the next section.

\subsection{Structural identifiability and observability of individual parameters and states}\label{ssec:SI_and_SO_individual}
We generalise Renardy et al.'s definition of structural identifiability~\cite{renardy2022structural} to additionally account for structural observability by explicitly including both parameters and states. Intuitively, a parameter is \emph{globally} identifiable if changes in that parameter necessarily lead to changes in the outputs for admissible states and parameters (see Defn.~\ref{def:global_SI}).

\begin{definition}[Global structural identifiability of a parameter]
\label{def:global_SI}
A parameter $\theta_{\ell}$ where $\ell\in\{1,\ldots,p\}$ is \emph{globally structurally identifiable} if any change in $\theta_{\ell}$ leads to a change in the outputs $\mathbf{y}(t,\mathbf{x},\boldsymbol{\theta})$ for all admissible states and parameters:
\begin{equation}
    \mathbf{y}(t,\mathbf{x},\boldsymbol{\theta}) = \mathbf{y}(t,\mathbf{x},\boldsymbol{\theta}^\star) \quad \forall t\in\mathbb{R}
    \quad \implies \quad \theta_{\ell} = \theta_{\ell}^{\star}.
\end{equation}
\end{definition}

Locally, a parameter is identifiable if small changes around a point in parameter space lead to changes in the outputs (see Defn.~\ref{def:local_SI}).

\begin{definition}[Local structural identifiability of a parameter]
\label{def:local_SI}
A parameter $\theta_{\ell}$ where $\ell\in\{1,\ldots,p\}$ is \emph{locally structurally identifiable} in $\mathcal{B}_{\boldsymbol{\theta}}(\boldsymbol{\theta}_C)$ if, for all $\boldsymbol{\theta},\boldsymbol{\theta}^\star \in \mathcal{B}_{\boldsymbol{\theta}}(\boldsymbol{\theta}_C)$ and admissible states $\mathbf{x}$,
\begin{equation}
    \mathbf{y}(t,\mathbf{x},\boldsymbol{\theta}) = \mathbf{y}(t,\mathbf{x},\boldsymbol{\theta}^\star) \quad \forall t\in\mathbb{R}
    \quad \implies \quad \theta_{\ell} = \theta_{\ell}^\star.
\end{equation}
\end{definition}

Similarly, a state is \emph{globally} observable if any change in the state leads to a change in outputs (see Defn.~\ref{def:global_SO}).

\begin{definition}[Global structural observability of a state]
\label{def:global_SO}
A state $x_i$ where $i\in\{1,\ldots,n\}$ is \emph{globally structurally observable} if any change in $x_i$ results in a change in the outputs for all admissible parameter values:
\begin{equation}
    \mathbf{y}(t,\mathbf{x},\boldsymbol{\theta}) = \mathbf{y}(t,\mathbf{x}^\star,\boldsymbol{\theta}) \quad \forall t\in\mathbb{R}
    \quad \implies \quad x_i = x_i^\star.
\end{equation}
\end{definition}

Locally, a state is observable if small deviations around a point in state space lead to changes in outputs (see Defn.~\ref{def:local_SO}).

\begin{definition}[Local structural observability of a state]
\label{def:local_SO}
A state $x_i$ where $i\in\{1,\ldots,n\}$ is \emph{locally structurally observable} in $\mathcal{B}_{\mathbf{x}}(\mathbf{x}_C)$ if, for all $\mathbf{x},\mathbf{x}^\star \in \mathcal{B}_{\mathbf{x}}(\mathbf{x}_C)$ and admissible parameters $\boldsymbol{\theta}$,
\begin{equation}
    \mathbf{y}(t,\mathbf{x},\boldsymbol{\theta}) = \mathbf{y}(t,\mathbf{x}^\star,\boldsymbol{\theta}) \quad \forall t\in\mathbb{R}
    \quad \implies \quad x_i = x_i^\star.
\end{equation}
\end{definition}

Even if individual parameters or states are not identifiable or observable on their own, specific combinations of them may be. We now turn to methods for characterising such identifiable and observable combinations.

\subsection{Structurally identifiable and observable combinations of parameters and states}
\label{ssec:SI_and_SO_comb}

We adapt Meshkat et al.'s notion of identifiable parameter combinations~\cite{meshkat2014global} to the framing of structural identifiability presented by Renardy et al.~\cite{renardy2022structural}, and generalise it to include observable combinations as well. A parameter combination is a function $\phi(\boldsymbol{\theta})$ of the parameters. Intuitively, a parameter combination is identifiable if changes in the combination necessarily lead to changes in the outputs for admissible states and parameters (see Defn.~\ref{def:global_param_comb}).

\begin{definition}[Global structural identifiability of a parameter combination]
\label{def:global_param_comb}
A function $\phi:\mathbb{R}^{p}\mapsto\mathbb{R}$ is \emph{globally structurally identifiable} if
\begin{equation}
    \mathbf{y}(t,\mathbf{x},\boldsymbol{\theta}) = \mathbf{y}(t,\mathbf{x},\boldsymbol{\theta}^\star) \quad \forall t\in\mathbb{R}
    \quad \implies \quad \phi(\boldsymbol{\theta}) = \phi(\boldsymbol{\theta}^\star).
\end{equation}
\end{definition}

Locally, a parameter combination is identifiable if small changes around a point in parameter space lead to changes in outputs (see Defn.~\ref{def:local_param_comb}).

\begin{definition}[Local structural identifiability of a parameter combination]
\label{def:local_param_comb}
A function $\phi:\mathbb{R}^{p}\mapsto\mathbb{R}$ is \emph{locally structurally identifiable} in $\mathcal{B}_{\boldsymbol{\theta}}(\boldsymbol{\theta}_C)$ if, for all $\boldsymbol{\theta},\boldsymbol{\theta}^\star \in \mathcal{B}_{\boldsymbol{\theta}}(\boldsymbol{\theta}_C)$,
\begin{equation}
    \mathbf{y}(t,\mathbf{x},\boldsymbol{\theta}) = \mathbf{y}(t,\mathbf{x},\boldsymbol{\theta}^\star) \quad \forall t\in\mathbb{R}
    \quad \implies \quad \phi(\boldsymbol{\theta}) = \phi(\boldsymbol{\theta}^\star).
\end{equation}
\end{definition}

A parameter-state combination is a function $\phi(\boldsymbol{\theta},\mathbf{x})$ of both parameters and states. Intuitively, a parameter-state combination is observable if changes in the combination necessarily lead to changes in outputs for admissible states and parameters (see Defn.~\ref{def:global_paramstate_comb}).

\begin{definition}[Global structural observability of a parameter-state combination]
\label{def:global_paramstate_comb}
A function $\phi:\mathbb{R}^{p}\times\mathbb{R}^{n}\mapsto\mathbb{R}$ is \emph{globally structurally observable} if
\begin{equation}
    \mathbf{y}(t,\mathbf{x},\boldsymbol{\theta}) = \mathbf{y}(t,\mathbf{x}^\star,\boldsymbol{\theta}^\star) \quad \forall t\in\mathbb{R}
    \quad \implies \quad \phi(\boldsymbol{\theta},\mathbf{x}) = \phi(\boldsymbol{\theta}^\star,\mathbf{x}^\star).
\end{equation}
\end{definition}

Locally, a parameter-state combination is observable if small deviations around a point in parameter-state space lead to changes in outputs (see Defn.~\ref{def:local_paramstate_comb}).

\begin{definition}[Local structural observability of a parameter-state combination]
\label{def:local_paramstate_comb}
A function $\phi:\mathbb{R}^{p}\times\mathbb{R}^{n}\mapsto\mathbb{R}$ is \emph{locally structurally observable} in $\mathcal{B}_{\boldsymbol{\theta}}(\boldsymbol{\theta}_C)\times \mathcal{B}_{\mathbf{x}}(\mathbf{x}_C)$ if, for all $(\boldsymbol{\theta},\mathbf{x}), (\boldsymbol{\theta}^\star,\mathbf{x}^\star) \in \mathcal{B}_{\boldsymbol{\theta}}(\boldsymbol{\theta}_C)\times \mathcal{B}_{\mathbf{x}}(\mathbf{x}_C)$,
\begin{equation}
    \mathbf{y}(t,\mathbf{x},\boldsymbol{\theta}) = \mathbf{y}(t,\mathbf{x}^\star,\boldsymbol{\theta}^\star) \quad \forall t\in\mathbb{R}
    \quad \implies \quad \phi(\boldsymbol{\theta},\mathbf{x}) = \phi(\boldsymbol{\theta}^\star,\mathbf{x}^\star).
\end{equation}
\label{def:obs_param_state_comb}
\end{definition}

As can be seen from these definitions, the conservation of observed outputs at each time point is key to understanding structural identifiability and observability. This conservation property provides a direct link to the concept of full symmetries. Building on this concept, our parameter–state symmetries (presented in Section~\ref{sec:results}) constitute a special case of full symmetries, which is defined at the end of this section. Full symmetries, in turn, are extensions of so-called classical symmetries, which we discuss next.

\subsection{Using classical Lie symmetries to uncover invariant model structures}\label{sec:classical_symmetries}

Lie symmetries constitute an immensely powerful tool widely used in theoretical physics, underpinning multiple Nobel Prizes, yet they are rarely utilised in fields such as mathematical biology. Classical Lie symmetries are particularly useful in studying differential equations. Their main applications are to find analytical solutions of ODEs (and PDEs), to reduce systems of ODEs to a simpler form, and, crucially, to construct classes of ODEs or PDEs that share common symmetries. 

For readers unacquainted with classical Lie theory, textbooks~\cite{hydon2000symmetry,bluman1989symmetries,olver2000applications} often present dense notation that can obscure the main ideas. The purpose of this section is to introduce these concepts in a clear and accessible manner, focusing on the intuition and fundamental properties of classical Lie symmetries before we extend them to full symmetries acting on parameters relevant to local structural identifiability and observability analyses.

\subsubsection{Classical Lie symmetries are transformations mapping solutions to solutions}\label{ssec:classical_symmetry}

A symmetry is a Lie transformation that maps solutions of ODEs to other solutions. To express this sentiment mathematically, we introduce a framework for accounting for derivatives, since ODEs are defined by derivatives. The required machinery involves prolongations, Jet spaces, and solution manifolds.

Consider an $n$-dimensional ODE system for states $\mathbf{x}\in\mathbb{R}^{n}$ depending on an independent variable $t\in\mathbb{R}$:
\begin{equation}
  \dfrac{\mathrm{d}\mathbf{x}}{\mathrm{d}t} = \mathbf{f}(t,\mathbf{x}).
  \label{eq:ODE_sys_classical}
\end{equation}
For now, we treat any parameters arising in these equations as constants meaning that the symmetries do not act on them. Let $\Gamma_{\varepsilon}$ be a one-parameter Lie transformation $\Gamma_\varepsilon:\mathbb{R}\times\mathbb{R}^{n}\to \mathbb{R}\times\mathbb{R}^{n}$:
\begin{equation}
  \Gamma_\varepsilon:(t,\mathbf{x})\mapsto (\hat{t}(t,\mathbf{x},\varepsilon), \hat{\mathbf{x}}(t,\mathbf{x},\varepsilon)),
  \label{eq:Lie_classical}
\end{equation}
with the trivial transformation obtained by setting $\varepsilon=0$:
\begin{equation}
\hat{t}(t,\mathbf{x},\varepsilon=0)=t,\quad \hat{\mathbf{x}}(t,\mathbf{x},\varepsilon=0)=\mathbf{x}.
\end{equation}
For simplicity, we write $\hat{t}(\varepsilon)$ instead of $\hat{t}(t,\mathbf{x},\varepsilon)$ and $\hat{\mathbf{x}}(\varepsilon)$ instead of $\hat{\mathbf{x}}(t,\mathbf{x},\varepsilon)$.

As we will see later, Lie transformations are symmetries of ODEs if they map solutions to other solutions. Since ODEs are defined by derivatives, Lie transformations must be extended to account for derivatives. To define the notion of a transformed derivative $\widehat{(\mathrm{d}\mathbf{x}/\mathrm{d}t)}(\varepsilon)$, we introduce the total derivative
\begin{equation}
  D_t = \partial_t + \sum_{i=1}^{n} \dfrac{\mathrm{d}x_{i}}{\mathrm{d}t} \partial_{x_i},
  \label{eq:total_derivative}
\end{equation}
so that the transformed derivative is given by
\begin{equation}
\widehat{\dfrac{\mathrm{d}\mathbf{x}}{\mathrm{d}t}}(\varepsilon) = \frac{D_t \hat{\mathbf{x}}(\varepsilon)}{D_t \hat{t}(\varepsilon)}.
\end{equation}
This, in turn, allows us to define the \textit{first prolongation} of the Lie symmetry $\Gamma_\varepsilon^{(1)}$:
\begin{equation}
\Gamma_\varepsilon^{(1)}:\left(t,\mathbf{x},\dfrac{\mathrm{d}\mathbf{x}}{\mathrm{d}t}\right)\mapsto \left(\hat{t}(\varepsilon), \hat{\mathbf{x}}(\varepsilon), \widehat{\dfrac{\mathrm{d}\mathbf{x}}{\mathrm{d}t}}(\varepsilon)\right),
\end{equation}
acting naturally on the \textit{first Jet space} $\mathcal{J}^{(1)}$
\begin{equation}
\mathcal{J}^{(1)} = \mathbb{R}\times\mathbb{R}^{n}\times\mathbb{R}^{n},
\end{equation}
where the \textit{solution manifold} $\mathcal{S}$ of the ODE system in Eq. \eqref{eq:ODE_sys_classical} is a subvariety of this Jet space according to
\begin{equation}
\mathcal{S} \coloneqq \left\{\left(t,\mathbf{x},\dfrac{\mathrm{d}\mathbf{x}}{\mathrm{d}t}\right)\in \mathcal{J}^{(1)} : \Delta =   \dfrac{\mathrm{d}\mathbf{x}}{\mathrm{d}t} - \mathbf{f}(t,\mathbf{x}) = \mathbf{0}\right\}.
\end{equation}
Each prolonged transformation $\Gamma_\varepsilon^{(1)}$ is uniquely defined by $\Gamma_\varepsilon$, so working with prolongations is equivalent to working with the original transformations.

Given this precise notion of prolongations, we can express the sentiment ``\textit{a symmetry maps solutions of ODEs to other solutions}'' mathematically. A Lie transformation $\Gamma_{\varepsilon}$ is a \textit{symmetry of the ODE system} in Eq. \eqref{eq:ODE_sys_classical} if its first prolongation preserves the solution manifold, i.e.,
\begin{equation}
  \Gamma^{(1)}_{\varepsilon} : \mathcal{S} \mapsto \mathcal{S},
\end{equation}
or equivalently,
\begin{equation}
  \widehat{\dfrac{\mathrm{d}\mathbf{x}}{\mathrm{d}t}}(\varepsilon) = \mathbf{f}(\hat{t}(\varepsilon), \hat{\mathbf{x}}(\varepsilon)) \quad \text{whenever} \quad \dfrac{\mathrm{d}\mathbf{x}}{\mathrm{d}t} = \mathbf{f}(t,\mathbf{x}).
  \label{eq:sym_cond}  
\end{equation}
In practice, however, we do not work directly with these \textit{symmetry conditions} in Eq. \eqref{eq:sym_cond} as we can formulate an alternative set of equations that are easier to solve and that completely characterise the symmetries. 

\subsubsection{Linearisation of Lie symmetries yields the key constituent of the analysis: infinitesimals}

A key property of Lie transformations, including symmetries, is that they can be Taylor expanded in the transformation parameter $\varepsilon$. This is guaranteed because Lie transformations are $\mathcal{C}^{\infty}$ diffeomorphisms. The linear terms in these expansions are called the \textit{infinitesimals} and are central to the study of differential equations using Lie symmetries, since each Lie transformation is completely characterised by its infinitesimals. Explicitly, we write
\begin{align}
\hat{t}(\varepsilon) &= t + \underbrace{\xi(t,\mathbf{x})}_{\left.\frac{d\hat{t}}{d\varepsilon}\right|_{\varepsilon=0}} \varepsilon + \mathcal{O}(\varepsilon^2),\\
\hat{x}_i(\varepsilon) &= x_i + \underbrace{\eta_i(t,\mathbf{x})}_{\left.\frac{d\hat{x}_i}{d\varepsilon}\right|_{\varepsilon=0}} \varepsilon + \mathcal{O}(\varepsilon^2),\quad i=1,\dots,n.
\end{align}
for the Taylor expansions around $\varepsilon=0$. Here, the functions $\xi(t,\mathbf{x})$ and $\eta_{i}(t,\mathbf{x})$ where $i=1,\ldots,n$ are referred to as the \textit{infinitesimals}. These infinitesimals \textit{generate the transformation} $\Gamma_{\varepsilon}:(t,\mathbf{x})\mapsto(\hat{t}(\varepsilon),\hat{\mathbf{x}}(\varepsilon))$ through the ODE system:
\begin{align}
\frac{d\hat{t}}{d\varepsilon} &= \xi(\hat{t},\hat{\mathbf{x}}), \quad \hat{t}(\varepsilon=0)=t,\\
\frac{d\hat{x}_i}{d\varepsilon} &= \eta_i(\hat{t},\hat{\mathbf{x}}), \quad \hat{x}_i(\varepsilon=0)=x_i,\quad i=1,\dots,n.
\end{align}
Thus, by finding the infinitesimals we obtain the corresponding pointwise Lie transformation $\Gamma_{\varepsilon}$ by solving the ODE system above. In other words, any Lie transformation is completely characterised by its infinitesimals.  A convenient notation based on these infinitesimals is the so-called \textit{infinitesimal generator of the Lie group} $X$:
\begin{equation}
  X = \xi \partial_t + \sum_{i=1}^{n} \eta_i \partial_{x_i}\,,
  \label{eq:X}
\end{equation}
and we say that \textit{$\Gamma_{\varepsilon}$ is generated by $X$}. Just like $\Gamma_\varepsilon$ has a first prolongation $\Gamma_\varepsilon^{(1)}$, the generator $X$ has a first prolongation $X^{(1)}$ given by
\begin{equation}
  X^{(1)} = X + \sum_{i=1}^{n} \eta_i^{(1)} \partial_{\mathrm{d}x_{i}/\mathrm{d}t},
  \label{eq:X_1}
\end{equation}
where the prolonged infinitesimals $\eta_i^{(1)}$ are given by the \textit{prolongation formula}:
\begin{equation}
  \eta_i^{(1)} = D_t \eta_i - \dfrac{\mathrm{d}x_{i}}{\mathrm{d}t}D_t \xi.
  \label{eq:prolongation_formula}
\end{equation}
The \textit{linearised symmetry conditions} are then expressed as
\begin{equation}
  X^{(1)}\left(  \dfrac{\mathrm{d}\mathbf{x}}{\mathrm{d}t} - \mathbf{f}(t,\mathbf{x})\right) = \mathbf{0} \quad \text{whenever} \quad   \dfrac{\mathrm{d}\mathbf{x}}{\mathrm{d}t} - \mathbf{f}(t,\mathbf{x}) = \mathbf{0},
  \label{eq:lin_sym}
\end{equation}
which correspond to the $\mathcal{O}(\varepsilon)$ terms of the Taylor expansions of the symmetry conditions in Eq. \eqref{eq:sym_cond} around $\varepsilon=0$. In essence, a symmetry analysis consists of solving these linearised conditions, which are PDEs, for the infinitesimals $\xi(t,\mathbf{x})$ and $\eta_{i}(t,\mathbf{x})$ where $i\in\left\{1,\ldots,n\right\}$. Recently, these classical symmetries have been extended to account for parameters which is relevant in the study of structural identifiability and observability. 

\subsection{Extending classical Lie symmetries to full symmetries connects symmetry and structural properties}\label{ssec:full_symmetry}
We are now at a point where we can connect Lie symmetries with local structural identifiability and observability. Structural identifiability (see Sections \ref{ssec:SI_and_SO_individual} and \ref{ssec:SI_and_SO_comb}) is related to changes in parameters that preserve observed outputs. In the context of identifiability, it is therefore natural to extend classical Lie symmetries to act on parameters \textit{in addition} to the independent (e.g. time $t$) and dependent (e.g. states $\mathbf{x}$) variables. To this end, we extend the classical symmetries $\Gamma_{\varepsilon}$ in Eq. \eqref{eq:Lie_classical} in the following way
\begin{equation}
  \Gamma_\varepsilon:(t,\mathbf{x},\boldsymbol{\theta})\mapsto (\hat{t}(\varepsilon), \hat{\mathbf{x}}(\varepsilon),\hat{\boldsymbol{\theta}}(\varepsilon)),
  \label{eq:full_symmetry}
\end{equation}
where we have additional transformed coordinates $\hat{\boldsymbol{\theta}}(\varepsilon)=\hat{\boldsymbol{\theta}}(t,\mathbf{x},\boldsymbol{\theta},\varepsilon)$ for the parameters. Infinitesimally speaking, the generating vector field becomes  
\begin{equation}
  X = \xi \partial_t + \sum_{i=1}^{n} \eta_i \partial_{x_i} + \sum_{\ell=1}^{p} \chi_{\ell} \partial_{\theta_\ell}\,,
  \label{eq:X_param}
\end{equation}
which contains an extra $p$ parameter infinitesimals $\chi_{\ell}(t,\mathbf{x},\boldsymbol{\theta})$ where $\ell\in\{1,\ldots,p\}$. Otherwise, the corresponding first prolongation $X^{(1)}$ is defined analogously to the classical counterpart in Eq. \eqref{eq:X_1} and so are the corresponding linearised symmetry conditions which are given by Eq. \eqref{eq:lin_sym}. However, there is a crucial conceptual difference between classical and full symmetries.

Whereas classical symmetries map solutions to other solutions, full symmetries do not satisfy this criteria as they transform parameters and thereby effectively change the model. In other words, full symmetries map solutions of one ODE, defined by a set of parameters, to a solution of another ODE, defined by another set of parameters. The linearised symmetry conditions of full symmetries ensure that the \textit{structure of a given class of models} is preserved, and thus these structure-preserving linearised symmetry conditions ensure that we map solutions to other solutions \textit{within a class of models}. It is the linearised symmetry conditions combined with an additional set of equations that enable the study of structural identifiability using full symmetries.

Previous work on full symmetries~\cite{massonis2020finding,villaverde2022symmetries,merkt2015higher}, have introduced an essential equation in the context of structural identifiability ensuring that \textit{observed outputs are preserved}. The mathematical equations for this condition are given by
\begin{equation}
  X(\mathbf{y}) = \mathbf{0}.
  \label{eq:output_conserved}
\end{equation}
These extra conditions express the fundamental requirement that the observed outputs remain invariant under transformations by full symmetries. Given the structure-preserving linearised symmetry conditions together with the output-preserving equations in Eq. \eqref{eq:output_conserved}, full symmetries have been used to analyse structural identifiability and observability. Nevertheless, the exact link between structural identifiability, observability, and the properties of full symmetries remains elusive, although recent work has begun to shed light on a potential answer to this challenge.

In other contexts, the explicit link between structural identifiability and observability, on the one hand, and full symmetries, on the other, is linked to the notion of invariance. This intuition was originally proposed by Castro and deBoer~\cite{castro2020structuralIdentifiability} which studied a particular class of full symmetries corresponding to scalings. They proposed that a parameter is structurally identifiable if the only possible value of the scaling factor of this parameter that leaves observed outputs unchanged is one. Similarly, they proposed that a state is structurally observable if the only possible value of the scaling factor of this state that leaves observed outputs unchanged is one. Essentially, these two statements express the idea that parameters and states are structurally identifiable and observable, respectively, if they are invariant under scalings. In our recent work~\cite{borgqvist2024framingglobalstructuralidentifiability}, we generalised and formalised the same notion of parameter identifiability by showing that the invariance of locally structurally identifiable parameter combinations holds for any parameter symmetry transforming parameters of the output system, not just scalings. Specifically, we showed that a parameter is structurally identifiable if and only if it is a so-called universal parameter invariant of all parameter symmetries of the output system. In fact, the identifiable parameter combinations of that model correspond to the universal parameter invariants. The parameter symmetries in our previous work were restricted to the output system, where states have been eliminated, and thus we could only analyse structural identifiability and not structural observability then. In this work, we generalise this idea further by studying the full model in Eq. \eqref{eq:ODE_sys_full} consisting of a state-dependent ODE system together with equations for observed ouputs, and we introduce a special type of full symmetry which we call \emph{parameter-state symmetries}.

The invariance of the observed outputs provides the key link between the symmetry framework introduced above and the notions of structural identifiability and structural observability. In the preceding discussion, we first introduced classical Lie symmetries acting on independent variables (e.g.\ time $t$) and dependent variables (e.g.\ states $\mathbf{x}$), and then extended this framework to full symmetries that also act on the parameters of the model. Imposing the condition $X(\mathbf{y})=\mathbf{0}$ on these full symmetries ensures that these transformations preserve the observed outputs, thereby identifying the class of transformations that are relevant for analysing structural identifiability and structural observability. The analysis of these structural properties is captured by \textit{parameter-state symmetries}, a special case of full symmetries that builds on the invariance of the observed outputs \textit{at each time point} and that will be introduced in the next section.

\section{Results}
\label{sec:results}
The results are presented in three parts, each building on the mathematical preliminaries introduced in Section~\ref{sec:math_preliminaries}. 
First, we formally define parameter-state symmetries and associated differential invariants, establishing the minimal set of transformations that preserve the observed outputs at each time point and are relevant for analysing structural identifiability and structural observability. Next, we state and prove our main theorems showing that the local structural properties of a model can be fully characterised by studying so-called \textit{universal invariants} of that model, providing a rigorous link between Lie symmetries and the classical notions of identifiability and observability. Finally, we analyse the local structural identifiability and local structural observability of four concrete models of biological systems by calculating the universal invariants of these models.

\subsection{Parameter-state symmetries as transformations preserving time and observed outputs}
\label{sec:defs_param_state_sym}

We define \textit{parameter-state symmetries} as the minimal class of continuous transformations of states and parameters that preserve the observed outputs at each time point (see Defn. \ref{def:SI_SO}). 

\begin{definition}[Parameter-state symmetries of the state-dependent system with observed outputs]
Consider a one-parameter family of transformations
\begin{equation}
  \Gamma_{\varepsilon}^{\mathbf{x},\boldsymbol{\theta}} : (t,\mathbf{x},\boldsymbol{\theta}) \mapsto 
  \big(t,\hat{\mathbf{x}}(t,\mathbf{x},\boldsymbol{\theta},\varepsilon), \hat{\boldsymbol{\theta}}(\boldsymbol{\theta},\varepsilon) \big).
\end{equation}
The generator of this transformation is
\begin{equation}
  X = \sum_{i=1}^{n} \eta_i(t,\mathbf{x},\boldsymbol{\theta}) \, \partial_{x_i} + \sum_{\ell=1}^{p} \chi_\ell(\boldsymbol{\theta}) \, \partial_{\theta_\ell}\,,
  \label{eq:X_param_state}
\end{equation}
where the time infinitesimal satisfies $\xi = 0$, ensuring that the transformation does not shift the independent variable $t$. This restriction guarantees that the observed outputs are preserved \emph{at each time point}, in accordance with the definitions of structural identifiability and structural observability, and prevents matching outputs at different times, which would violate these structural definitions.
The parameter infinitesimals $\chi_\ell$ for $\ell=1,\ldots,p$ depend solely on the parameters $\boldsymbol{\theta}$. The corresponding first prolongation of the generator is
\begin{equation}
  X^{(1)} = X + \sum_{i=1}^{n} \eta_i^{(1)}\left(t,\mathbf{x},\dfrac{\mathrm{d}\mathbf{x}}{\mathrm{d}t},\boldsymbol{\theta}\right) \, \partial_{\mathrm{d}x_{i}/\mathrm{d}t},
\end{equation}
where $\eta_i^{(1)}$ is defined by the prolongation formula in Eq. \eqref{eq:prolongation_formula} in Section~\ref{sec:math_preliminaries}.

We say that $\Gamma_{\varepsilon}^{\mathbf{x},\boldsymbol{\theta}}$ is a \textbf{parameter-state symmetry} of the state-dependent ODE system together with observed output equations in Eq.~\eqref{eq:ODE_sys_full} if it satisfies both the \textbf{linearised symmetry conditions}
\begin{equation}
  X^{(1)}\left(\dfrac{\mathrm{d}\mathbf{x}}{\mathrm{d}t} - \mathbf{f}(t,\mathbf{x},\boldsymbol{\theta})\right) = \mathbf{0} \quad \text{whenever } \dfrac{\mathrm{d}\mathbf{x}}{\mathrm{d}t} = \mathbf{f}(t,\mathbf{x},\boldsymbol{\theta})\,,
  \label{eq:lin_sym_param_state}
\end{equation}
together with the \textbf{output invariance conditions}
\begin{equation}
  X(\mathbf{y}) = \mathbf{0}\,.
  \label{eq:X_y_zero}
\end{equation}
\label{def:SI_SO}
\end{definition}

The \textit{linearised symmetry conditions}, the \textit{output invariance conditions} and the restrictions on the infinitesimals of the parameter-state symmetries together define the minimal class of continuous transformations of states and parameters that preserve the observed outputs at each time point. Importantly, parameter symmetries (Def. \ref{def:SI_SO}) constitute a special case of the previously defined full symmetries (Section \ref{ssec:full_symmetry}). In the full symmetry framework, no restrictions on the time or parameter infinitesimals are assumed, i.e., $\xi=\xi(t,\mathbf{x},\boldsymbol{\theta})$ and $\chi_\ell = \chi_\ell(t,\mathbf{x},\boldsymbol{\theta})$ for $\ell=1,\ldots,p$ as in Eq.~\eqref{eq:X_param}, and a full symmetry defined by a non-constant time infinitesimal will shift time. Thus, the output invariance conditions defining full symmetries, in general, preserve the observed ouputs while also allowing for time shifts. The key restriction of parameter-state symmetries is ultimately $\xi=0$ as it prohibits time shifts and together with the output invariance conditions it ensures that the oberved outputs are preserved \textit{at each time point}, with the first prolongation $X^{(1)}$ and linearised symmetry conditions defined analogously to the classical counterparts in Eq.~\eqref{eq:X_1} and Eq.~\eqref{eq:lin_sym}. This clarifies the previously elusive link between full symmetries and the local structural properties of the model, through the conservation of observed outputs at each time point, which defines parameter–state symmetries.

Following the definition of parameter-state symmetries, we now introduce \textit{differential invariants}, which are quantities preserved under the action of these parameter-state symmetries. Differential invariants (Defn. \ref{def:diff_inv}) can be used to analyse the local structural identifiability and structural observability of a model, as they remain unchanged under transformations by parameter-state symmetries.

\begin{definition}[Differential invariant of a parameter-state symmetry]
Let $\Gamma_{\varepsilon}^{\mathbf{x},\boldsymbol{\theta}}$ be a one-parameter family of parameter-state symmetries with generator $X$ as in Defn.~\ref{def:SI_SO}. A smooth function $I(t,\mathbf{x},\mathrm{d}\mathbf{x}/\mathrm{d}t,\boldsymbol{\theta})$ is called a \textbf{differential invariant} of $\Gamma_{\varepsilon}^{\mathbf{x},\boldsymbol{\theta}}$ if it satisfies
\begin{equation}
  X^{(1)}(I) = 0\,,
  \label{eq:differential_invariant}
\end{equation}
where $X^{(1)}$ denotes the first prolongation of $X$.
\label{def:diff_inv}
\end{definition}

By construction, differential invariants are constant under transformations by parameter-state symmetries. Consequently, they provide quantities that remain unchanged under the minimal class of transformations preserving the observed outputs at each time point. As we will show in the following theorems, these invariants can be used to fully characterise both the local structural identifiability and local structural observability of any model.

Depending on which structural property we wish to analyse, it is useful to categorise differential invariants according to their dependence on states and parameters. In the context of local structural identifiability, \textit{parameter invariants} are invariants that depend solely on the parameters $\boldsymbol{\theta}$, i.e., $I(\boldsymbol{\theta})$. In the context of local structural observability, \textit{state invariants} are invariants that depend solely on the states $\mathbf{x}$, i.e., $I(\mathbf{x})$. We can also define \textit{parameter-state invariants}, which depend on both the states $\mathbf{x}$ and parameters $\boldsymbol{\theta}$, i.e., $I(\mathbf{x},\boldsymbol{\theta})$. Moreover, some differential invariants are specific to particular parameter-state symmetries, while others are shared by \textit{all parameter-state symmetries} of the model. We refer to the latter type of invariant as a \textit{universal invariant} (Defn.~\ref{def:uni_inv}).

\begin{definition}[Universal invariants of the state-dependent system with observed outputs]
  Consider a model as in Eq.~\eqref{eq:ODE_sys_full} consisting of a state-dependent ODE system with equations for observed outputs . A smooth function $I(t,\mathbf{x},\mathrm{d}\mathbf{x}/\mathrm{d}t,\boldsymbol{\theta})$ that is a differential invariant (Defn.~\ref{def:diff_inv}) of \textbf{all} parameter-state symmetries of this model, i.e., $X^{(1)}(I) = 0$ for every generator $X$, is called a \textbf{universal invariant} of the model. 
\label{def:uni_inv}
\end{definition}

In analogy with the classification of differential invariants introduced above, we can distinguish between different types of universal invariants. In particular, \textit{universal parameter invariants} depend solely on the parameters $\boldsymbol{\theta}$, \textit{universal state invariants} depend solely on the states $\mathbf{x}$, and \textit{universal parameter-state invariants} depend on both the states and the parameters, i.e., $\mathbf{x}$ and $\boldsymbol{\theta}$. These categories of universal invariants will be central in the following section, where we formally link them to the local structural identifiability and local structural observability of a model.

\subsection{Using universal invariants to characterise local structural identifiability and observability of models}
\label{sec:theorems_SI_SO}
We now state our main theoretical results, which establish a rigorous connection between parameter-state symmetries, on the one hand, and the local structural identifiability and local structural observability of model, on the other. In particular, we show that the local structural properties of \textit{the model} — by which we mean the state-dependent ODE system together with its observed outputs as in Eq.~\eqref{eq:ODE_sys_full} — can be completely characterised by its universal invariants. Universal parameter invariants capture all structurally identifiable combinations of parameters, universal state invariants capture all structurally observable combinations of states, and the universal parameter-state invariants correspond to a new type of quantity which we refer to as a \textit{structurally observable combination}. These results formalise and generalise the original intuition of Castro and deBoer~\cite{castro2020structuralIdentifiability} that structurally identifiable parameters are invariant under scalings and that structurally observable states are also invariant under scalings. They also generalise our previous results showing that all structurally identifiable parameter combinations are universal parameter invariants of the parameter symmetries of the output system~\cite{borgqvist2024framingglobalstructuralidentifiability}, which is the higher-order output-equivalent system of the model where all states have been eliminated. In total, these results allow us to simultaneously analyse local structural identifiability and local structural observability through the universal invariants of the model.

For clarity, we note that these universal invariants can be classified according to their dependence on states and parameters, which directly determines whether they correspond to structurally identifiable or structurally observable quantities. Universal invariants depending solely on parameters correspond to structurally identifiable combinations, whereas universal invariants that depend on states, or on states together with parameters, correspond to structurally observable combinations. For instance, a quantity such as $\beta x$, involving both a parameter $\beta$ and a state $x$, is classified as a structurally observable combination (Defn. \ref{def:obs_param_state_comb} in Section \ref{ssec:SI_and_SO_comb}). Likewise, structurally identifiable quantities solely depending on parameters and structurally observable quantities which depend on states are distinguished from one another. This distinction ensures that the subsequent theorems can be interpreted consistently with the definitions of structural identifiability and structural observability (see Sections \ref{ssec:SI_and_SO_individual} and \ref{ssec:SI_and_SO_comb}).

We now show that the universal invariants of parameter-state symmetries fully characterise which parameters are locally structurally identifiable and which states are locally structurally observable. This establishes a rigorous connection between the parameter-state symmetries of any model on the form presented in Eq. \eqref{eq:ODE_sys_full} and its local structural properties. We begin by demonstrating that a parameter is locally structurally identifiable if and only if it is a universal parameter invariant (Theo. \ref{thm:SI_symmetries}).

\begin{theo}[Local structural identifiability in terms of universal parameter invariants]
    Let $\mathbf{y}(t,\mathbf{x},\boldsymbol{\theta})$ be a curve for the observed outputs of the model in Eq.~\eqref{eq:ODE_sys_full} where the states solve $\mathrm{d}\mathbf{x}/\mathrm{d}t=\mathbf{f}(t,\mathbf{x},\boldsymbol{\theta})$. A parameter $\theta_{\ell}\in\boldsymbol{\theta}$, ${\ell}\in\{1,\ldots,p\}$, is locally structurally identifiable if and only if it is a universal parameter invariant of the model.
\label{thm:SI_symmetries}
\end{theo}
\begin{proof}
\textbf{``$\Longrightarrow$''} Assume $\theta_{\ell}$ is locally structurally identifiable. We need to show that $\theta_\ell$ also is a universal parameter invariant. Take a radius $\rho_{\boldsymbol{\theta}}>0$ and a centre $\boldsymbol{\theta}_{C}\in\mathbb{R}^{p}$ defining the ball $\mathcal{B}_{\boldsymbol{\theta}}(\boldsymbol{\theta}_C)$ in parameter space where $\theta_\ell$ is identifiable. Let $\rho_{\mathbf{x}}>0$ and $\mathbf{x}_C\in\mathbb{R}^{n}$ define the corresponding ball $\mathcal{B}_{\mathbf{x}}(\mathbf{x}_C)$ in state space. By the definition of local identifiability (see Defn. \ref{def:local_SI} in Section \ref{ssec:SI_and_SO_individual}), we have
\begin{equation}
\mathbf{y}(t,\mathbf{x},\boldsymbol{\theta}) = \mathbf{y}(t,\mathbf{x}^\star,\boldsymbol{\theta}^\star) \implies \theta_\ell = \theta_\ell^\star
\quad \forall \boldsymbol{\theta},\boldsymbol{\theta}^\star \in \mathcal{B}_{\boldsymbol{\theta}}(\boldsymbol{\theta}_C), \forall \mathbf{x},\mathbf{x}^\star \in \mathcal{B}_{\mathbf{x}}(\mathbf{x}_C).
\label{eq:SI_condition}
\end{equation}
Based on Defn. \ref{def:SI_SO}, take a parameter-state symmetry $\Gamma^{\mathbf{x},\boldsymbol{\theta}}_\varepsilon:(t,\mathbf{x},\boldsymbol{\theta})\mapsto(t,\hat{\mathbf{x}}(\varepsilon),\hat{\boldsymbol{\theta}}(\varepsilon))$. Since this parameter-state symmetry satisfies the output-invariance condition in Eq. \eqref{eq:X_y_zero}, we have that $\mathbf{y}(t,\mathbf{x},\boldsymbol{\theta})=\mathbf{y}(t,\hat{\mathbf{x}}(\varepsilon),\hat{\boldsymbol{\theta}}(\varepsilon))$. Moreover, our parameter symmetry acts infinitesimally on all parameters in parameter space and all states in state space, i.e. on all $\boldsymbol{\theta}\in\mathbb{R}^{p}$ and all $\mathbf{x}\in\mathbb{R}^{n}$, in general, as well as on all local parameters and local states, i.e. on all $\boldsymbol{\theta}\in\mathcal{B}_{\boldsymbol{\theta}}(\boldsymbol{\theta}_C)$ and all $\mathbf{x}\in \mathcal{B}_{\mathbf{x}}(\mathbf{x}_C)$, in particular. Therefore, we can always find an interval $[\varepsilon_1,\varepsilon_2]$ for the transformation parameter $\varepsilon$ such that the transformed parameters $\hat{\boldsymbol{\theta}}(\varepsilon)$ remain in $\mathcal{B}_{\boldsymbol{\theta}}(\boldsymbol{\theta}_C)$ and the transformed states $\hat{\mathbf{x}}(\varepsilon)$ remain in $\mathcal{B}_{\mathbf{x}}(\mathbf{x}_C)$. Consequently, we have that
\begin{equation}
  \mathbf{y}(t,\mathbf{x},\boldsymbol{\theta}) = \mathbf{y}(t,\hat{\mathbf{x}}(\varepsilon),\hat{\boldsymbol{\theta}}(\varepsilon)) \implies \theta_\ell = \hat{\theta}_\ell(\varepsilon) \quad \forall \varepsilon \in [\varepsilon_1,\varepsilon_2]\,,
  \label{eq:symmetry_SI}
\end{equation}
and thus $\hat{\theta}_{\ell}(\varepsilon)$ is constant on the interval $[\varepsilon_1,\varepsilon_2]$. Moreover, since $\hat{\theta}_\ell(\varepsilon)$ is analytic, it must therefore be constant everywhere. Hence, $\theta_\ell$ is a parameter invariant of the parameter-state symmetry $\Gamma^{\mathbf{x},\boldsymbol{\theta}}_\varepsilon$. The same arguments hold for all parameter-state symmetries, and thus $\theta_\ell$ is a universal parameter invariant.\\  
\textbf{``$\Longleftarrow$''} Conversely, assume $\theta_\ell$ is a universal parameter invariant. We need to show that $\theta_{\ell}$ also is locally structurally identifiable. Specifically, we have that $\theta_{\ell}$ satisfies Eq. \eqref{eq:symmetry_SI} for all parameter-state symmetries as well as for all $\varepsilon\in\mathbb{R}$ and we need to show that we can define balls $\mathcal{B}_{\boldsymbol{\theta}}(\boldsymbol{\theta}_C)$ and $\mathcal{B}_{\mathbf{x}}(\mathbf{x}_C)$ in parameter space and state space, respectively, such that Eq. \eqref{eq:SI_condition} is satisfied within these balls. In other words, we need to show that there exists balls in parameter space and state space that all parameter-state symmetries, characterised by Eq. \eqref{eq:symmetry_SI}, cover. By definition, the output invariance conditions in Eq. \eqref{eq:X_y_zero} hold for all parameter-state symmetries of the model. Infinitesimally speaking, this implies that the entire set of parameter infinitesimals $\chi_{\ell}(\boldsymbol{\theta})$ where $\ell=1,\ldots,p$ derived from all parameter-state symmetries span all directions in parameter space for which the outputs are preserved. Since we can make the radii $\rho_{\boldsymbol{\theta}}>0$ and $\rho_{\mathbf{x}}>0$ in parameter space and state space, respectively, arbitrarily small, we can always define corresponding balls which the parameter-state symmetries cover, and hence $\theta_{\ell}$ is locally structurally identifiable.\\
\end{proof}

Building on this result, we can characterise \textit{structurally identifiable parameter combinations} (Defn. \ref{def:local_param_comb} in Section \ref{ssec:SI_and_SO_comb}) as universal parameter invariants as well. To this end, let $\mathbf{y}(t,\mathbf{x},\boldsymbol{\theta})$ be a curve for the observed outputs of a model of interest, and assume that all of its universal parameter invariants are collected in a vector $\boldsymbol{I}_{\boldsymbol{\theta}}\in\mathbb{R}^{\tilde{p}}$ for some $0<\tilde{p}<p$. Then, we can simply re-parametrise the observed outputs according to $\mathbf{y}(t,\mathbf{x},\boldsymbol{\theta})\mapsto\mathbf{y}(t,\mathbf{x},\boldsymbol{I}_{\boldsymbol{\theta}})$ and then apply Theo. \ref{thm:SI_symmetries} to the re-parametrised outputs. Consequently, all structurally identifiable parameter combinations correspond to the universal parameter invariants of the model at hand (Cor. \ref{cor:SI_symmetries}). 

\begin{cor}[Local structural identifiability in terms of universal parameter invariants]
The locally structurally identifiable parameter combinations of the model in Eq.~\eqref{eq:ODE_sys_full} are given by its universal parameter invariants.
\label{cor:SI_symmetries}
\end{cor}
The next result establishes the analogous statement for the local structural observability of states. More precisely, the locally structurally observable states are given by the universal state invariants of the model of interest (Theo. \ref{thm:SO_symmetries}). The proof for this result regarding structural observability follows the same reasoning as the case of structural identifiability of parameters presented in Theo. \ref{thm:SI_symmetries}, as we are essentially replacing parameters with states and universal parameter invariants with universal state invariants. For brevity, we therefore present a condensed version of the proof of Theo. \ref{thm:SO_symmetries}.

\begin{theo}[Local structural observability in terms of universal state invariants]
Let $\mathbf{y}(t,\mathbf{x},\boldsymbol{\theta})$ be a curve for the observed outputs of the model in Eq.~\eqref{eq:ODE_sys_full} where the states solve $\mathrm{d}\mathbf{x}/\mathrm{d}t=\mathbf{f}(t,\mathbf{x},\boldsymbol{\theta})$. A state $x_{i}\in\mathbf{x}$, $i\in\{1,\ldots,n\}$, is locally structurally observable if and only if it is a universal state invariant of the model.
\label{thm:SO_symmetries}
\end{theo}

\begin{proof}
\textbf{``$\Longrightarrow$''} Assume $x_{i}$ is locally structurally observable. Then there exist radii $\rho_{\mathbf{x}}>0$ and $\rho_{\boldsymbol{\theta}}>0$ with corresponding balls $\mathcal{B}_{\mathbf{x}}(\mathbf{x}_C)$ and $\mathcal{B}_{\boldsymbol{\theta}}(\boldsymbol{\theta}_C)$ such that
\begin{equation}
\mathbf{y}(t,\mathbf{x},\boldsymbol{\theta})=\mathbf{y}(t,\mathbf{x}^\star,\boldsymbol{\theta}^\star)
\;\Rightarrow\;
x_i=x_i^\star
\end{equation}
for all $(\mathbf{x},\boldsymbol{\theta})$ and $(\mathbf{x}^\star,\boldsymbol{\theta}^\star)$ in these neighbourhoods (see Defn. \ref{def:local_SO}).  

Consider now a parameter-state symmetry $\Gamma^{\mathbf{x},\boldsymbol{\theta}}_\varepsilon:(t,\mathbf{x},\boldsymbol{\theta})\mapsto(t,\hat{\mathbf{x}}(\varepsilon),\hat{\boldsymbol{\theta}}(\varepsilon))$. Since this transformation preserves the outputs, we have
\begin{equation}
\mathbf{y}(t,\mathbf{x},\boldsymbol{\theta})
=
\mathbf{y}(t,\hat{\mathbf{x}}(\varepsilon),\hat{\boldsymbol{\theta}}(\varepsilon)).
\end{equation}
As in the proof of (Theo. \ref{thm:SI_symmetries}), the infinitesimal action of the symmetry ensures that there exists an interval $[\varepsilon_1,\varepsilon_2]$ such that the transformed states and parameters remain within $\mathcal{B}_{\mathbf{x}}(\mathbf{x}_C)$ and $\mathcal{B}_{\boldsymbol{\theta}}(\boldsymbol{\theta}_C)$. Hence,
\begin{equation}
\mathbf{y}(t,\mathbf{x},\boldsymbol{\theta})
=
\mathbf{y}(t,\hat{\mathbf{x}}(\varepsilon),\hat{\boldsymbol{\theta}}(\varepsilon))
\;\Rightarrow\;
x_i=\hat{x}_i(\varepsilon)
\quad \forall \varepsilon\in[\varepsilon_1,\varepsilon_2].
\end{equation}
Thus $\hat{x}_i(\varepsilon)$ is constant on this interval, and by analyticity it is constant everywhere. Hence $x_i$ is a state invariant of the symmetry. Since the same argument holds for all parameter-state symmetries, $x_i$ is a universal state invariant.

\textbf{``$\Longleftarrow$''} Conversely, assume $x_i$ is a universal state invariant. Then $x_i$ is preserved under transformations by all parameter-state symmetries satisfying the output invariance conditions in Eq.~\eqref{eq:X_y_zero}. Infinitesimally, the corresponding state infinitesimals $\eta_i(\mathbf{x})$ describe all directions in state space that leave the outputs unchanged. By choosing the radii $\rho_{\mathbf{x}}$ and $\rho_{\boldsymbol{\theta}}$ arbitrarily small, we can define neighbourhoods $\mathcal{B}_{\mathbf{x}}(\mathbf{x}_C)$ and $\mathcal{B}_{\boldsymbol{\theta}}(\boldsymbol{\theta}_C)$ that are covered by the transformations of all parameter-state symmetries. Consequently, invariance of the outputs within these neighbourhoods implies invariance of $x_i$, and thus $x_i$ is locally structurally observable.
\end{proof}

Just as in the case of locally structural identifiability of parameter combinations, the reasoning in Theo. \ref{thm:SO_symmetries} can be extended to account for \textit{locally structurally observable parameter-state combinations} (Defn. \ref{def:obs_param_state_comb} in Section \ref{ssec:SI_and_SO_comb}). The corresponding result says that the \textit{structurally observable parameter-state combinations are given by the universal parameter-state invariants of the model of interest} (Cor. \ref{cor:SO_symmetries}). 

\begin{cor}[Local structural observability in terms of universal invariants]
The locally structurally observable parameter-state combinations $\phi(\boldsymbol{\theta},\mathbf{x})$ of the model in Eq.~\eqref{eq:ODE_sys_full} are given by its universal parameter-state invariants.
\label{cor:SO_symmetries}
\end{cor} 
The results above show that the local structural properties of a model can be characterised entirely through its universal invariants. In practice, analysing local structural identifiability and structural observability therefore reduces to calculating these universal invariants. Once the universal invariants of a model have been identified, all locally structurally identifiable parameter combinations and locally structurally observable state or parameter–state combinations follow directly. In the following section, we illustrate this procedure using several examples of mechanistic models of different biological systems.

\subsection{Examples of analyses of local structural identifiability and observability using parameter-state symmetries}
Here, we illustrate how our main theorems, which link local structural properties of models defined by observed outputs to universal invariants, are applied to actual mechanistic models. Specifically, we conduct a combined local structural identifiability and local structural observability analysis of four specific models on the same form as in Eq.~\eqref{eq:ODE_sys_full}. These four models, which increase in complexity, are a decoupled decay model, a linear model, a non-autonomous glucose-insulin model, and an epidemiological SEI model of tuberculosis spread. Importantly, the local structural identifiability of the same exact models were studied in our previous work~\cite{borgqvist2024framingglobalstructuralidentifiability} by means of parameter symmetries restricted to the corresponding output systems. The results presented here based on parameter-state symmetries of the full model generalise that analysis by incorporating the states into the symmetry framework. Interestingly enough, there is an equivalence between the parameter symmetries restricted to the output system analysed in our previous work and the parameter-state symmetries considered here. Just as the full model in Eq.~\eqref{eq:ODE_sys_full} is output equivalent to the higher-order output system, there is a similar equivalence between the parameter symmetries acting on the output system analysed in~\cite{borgqvist2024framingglobalstructuralidentifiability} and the parameter-state symmetries in this work, as the universal parameter invariants in both contexts are identical. Consequently, the locally structurally identifiable parameter combinations that are presented here recover exactly the results from~\cite{borgqvist2024framingglobalstructuralidentifiability}, while the locally structurally observable state or parameter-state combinations represent novel findings.

All calculations underlying the results presented here are detailed in the Appendices. In brief, the method for obtaining the parameter-state symmetries is straightforward. First, conditions on some of the infinitesimals are obtained by using the output invariance conditions in Eq. \eqref{eq:X_y_zero}. Thereafter, these infinitesimals are substituted into the linearised symmetry conditions in Eq. \eqref{eq:lin_sym_param_state} which are solved for the remaining infinitesimals. Lastly, the universal invariants $I$ are obtained by solving the equations $X(I)=0$, which are linear PDEs, for all obtained generators using the method of characteristics.

For the glucose-insulin model, we note that Massonis and Villaverde~\cite{massonis2020finding} also analysed the local structural properties of this model using full symmetries. To characterise the full symmetries having no a priori restrictions on the infinitesimals, these authors used an ansatz-based polynomial approach for finding the unknown infinitesimals. Their approach recovers the same fundamental scaling symmetry that we obtain here analytically from the linearised symmetry conditions, providing an independent confirmation of the structure of the parameter-state symmetries we find for this model. 

All calculations for the SEI model were conducted using the open-source symbolic solver SymPy~\cite{meurer2017sympy}, and the relevant scripts are available in the associated GitHub repository at \url{https://github.com/JohannesBorgqvist/parameter_symmetries_identifiability_and_observability/}.

\subsubsection{Decoupled decay model}

Consider a simple system consisting of two decoupled decay processes. Let \(u(t)\) and \(v(t)\) denote the concentrations of two substances at time \(t\), governed by

\begin{equation}
\begin{split}
    \frac{\mathrm{d}u}{\mathrm{d}t} &= \kappa_{1} - \lambda u\,,\\
    \frac{\mathrm{d}v}{\mathrm{d}t} &= \kappa_{2} - \lambda v\,,
\end{split}
\label{eq:decay_sys}
\end{equation}
with an observed output defined as

\begin{equation}
    y = u + v\,.
\label{eq:output_decay}
\end{equation}
The parameter vector is \(\boldsymbol{\theta} = (\kappa_{1}, \kappa_{2}, \lambda)^\top\), where \(\lambda\) is the common decay rate and \(\kappa_{1},\kappa_{2}\) are constant formation rates. We calculate the parameter-state symmetries of this model next and use them to analyse the local structural identifiability and structural observability of this model.

The family of infinitesimal generators is
\begin{equation}
    X = \eta_{u}(u,v,\boldsymbol{\theta}) (\partial_{u} - \partial_{v}) + \chi_{\kappa_{1}}(\boldsymbol{\theta}) (\partial_{\kappa_{1}} - \partial_{\kappa_{2}})\,,
\label{eq:X_decay_final}
\end{equation}
where \(\chi_{\kappa_{1}}:\mathbb{R}^{3} \to \mathbb{R}\) is an arbitrary function of the parameters, and where the state infinitesimal \(\eta_{u}:\mathbb{R}\times\mathbb{R}\times\mathbb{R}^{3}\mapsto\mathbb{R}\) solves the following PDE:
\begin{equation}
    \chi_{\kappa_{1}} = D_{t} \eta_{u} - \lambda \eta_{u}\,,
\label{eq:eta_u_decay_final}
\end{equation}
where $D_{t}$ is the total derivative introduced in Eq. \eqref{eq:total_derivative} in Section \ref{ssec:classical_symmetry}. The three \textit{universal invariants} are
\begin{equation}
    I_{1} = \lambda\,,\quad I_{2} = \kappa_{1} + \kappa_{2}\,,\quad I_{3} = u + v\,.
\label{eq:uni_inv_decay}
\end{equation}
From these invariants, the decay rate \(\lambda\) is \textit{locally structurally identifiable}, while the formation rates \(\kappa_{1}\) and \(\kappa_{2}\) are individually unidentifiable but their sum is \textit{locally structurally identifiable}. Similarly, \(u\) and \(v\) are unobservable individually, whereas the sum \(u+v\) is \textit{locally structurally observable} and constitutes a locally observable state combination. Importantly, this locally structurally observable state combination for the decoupled decay model corresponds exactly to the observed output, \(u+v\), highlighting that the symmetry analysis recovers the parameter-state combinations that are directly observed in the outputs. Having calculated the universal invariants, we study particular transformations of a parameter-state symmetry generated by a specific choice of generator. 

For the specific choice \(\eta_{u} = 1/\lambda\) and \(\chi_{\kappa_{1}} = 1\) of the infinitesimals, the corresponding vector field becomes

\begin{equation}
    X = \frac{1}{\lambda} (\partial_{u} - \partial_{v}) + (\partial_{\kappa_{1}} - \partial_{\kappa_{2}})\,,
\label{eq:X_decay_particular}
\end{equation}
which generates the following parameter-state symmetry
\begin{equation}
    \Gamma^{\mathbf{x},\boldsymbol{\theta}}_{\varepsilon}:(t,u,v,\kappa_{1},\kappa_{2},\lambda) \mapsto \left(t, u + \frac{\varepsilon}{\lambda}, v - \frac{\varepsilon}{\lambda}, \kappa_{1} + \varepsilon, \kappa_{2} - \varepsilon, \lambda\right)\,.
\label{eq:parameter_symmetry_sum}
\end{equation}
Transformations by \(\Gamma^{\mathbf{x},\boldsymbol{\theta}}_{\varepsilon}\) correspond to translations of the individual states and formation rates with opposite signs, while \textit{preserving the model structure} in Eq.~\eqref{eq:decay_sys}. As expected, transformations by the parameter-state symmetry $\Gamma^{\mathbf{x},\boldsymbol{\theta}}_{\varepsilon}$ in Eq.~\eqref{eq:parameter_symmetry_sum} preserves the structure of the decay model in Eq.~\eqref{eq:decay_sys}, which is illustrated by the following calculations:
\begin{equation}
    \begin{split}
        \widehat{\dfrac{\mathrm{d}u}{\mathrm{d}t}}(\varepsilon)&=\hat{\kappa_{1}}(\varepsilon)-\hat{\lambda}(\varepsilon){\hat{u}(\varepsilon)}\\
        \widehat{\dfrac{\mathrm{d}v}{\mathrm{d}t}}(\varepsilon)&=\hat{\kappa_{2}}(\varepsilon)-\hat{\lambda}(\varepsilon){\hat{v}(\varepsilon)}\\        
    \end{split}\quad\Longrightarrow\quad\begin{split}
        \dfrac{\mathrm{d}u}{\mathrm{d}t}&=(\kappa_{1}+\varepsilon)-\lambda{\left(u+\dfrac{\varepsilon}{\lambda}\right)}\\
        \dfrac{\mathrm{d}v}{\mathrm{d}t}&=(\kappa_{2}-\varepsilon)-\lambda{\left(v-\dfrac{\varepsilon}{\lambda}\right)}\\        
    \end{split}\quad\Longrightarrow\quad
        \begin{split}
        \dfrac{\mathrm{d}u}{\mathrm{d}t}&=\kappa_{1}-\lambda{u}\\
        \dfrac{\mathrm{d}v}{\mathrm{d}t}&=\kappa_{2}-\lambda{v}\\        
    \end{split}\,.
    \label{eq:decay_sys_preserved}
\end{equation}
In essence, Eq. \eqref{eq:decay_sys_preserved} demonstrates that the generated parameter-state symmetry satisfies the symmetry conditions in Eq. \eqref{eq:sym_cond} in Section \ref{ssec:classical_symmetry}. This example illustrates how parameter-state symmetries simultaneously reveal locally structurally identifiable parameter combinations and locally structurally observable state combinations.

\subsubsection{A linear model}

We now consider a simple linear system consisting of two coupled states, \(x(t)\) and \(z(t)\), whose dynamics are governed by
\begin{equation}
    \begin{split}
        \frac{\mathrm{d}x}{\mathrm{d}t} &= a x + b z\,,\\
        \frac{\mathrm{d}z}{\mathrm{d}t} &= c z\,,\\
    \end{split}
    \label{eq:linear_sys}
\end{equation}
with an observed output given by
\begin{equation}
    y = x\,,
    \label{eq:output_linear}
\end{equation}
and the parameter vector \(\boldsymbol{\theta} = (a, b, c)^\top\). The family of generating vector fields for this model is
\begin{equation}
    X = \chi_{b}(\boldsymbol{\theta}) \left( b \partial_{b} - z \partial_{z} \right)\,,
    \label{eq:X_linear_final}
\end{equation}
where \(\chi_{b}:\mathbb{R}^{3} \to \mathbb{R}\) is an arbitrary function of the parameters. This vector field generates the full set of parameter-state symmetries for the linear model.

The universal invariants of the linear model are
\begin{equation}
    I_1 = a\,, \quad I_2 = c\,, \quad I_3 = x\,, \quad I_4 = b z\,.
    \label{eq:uni_inv_linear}
\end{equation}
From these invariants, we see that the parameters \(a\) and \(c\) are \textit{locally structurally identifiable}, while the state \(x\), corresponding to the observed output, is \textit{locally structurally observable}. This is in agreement with the previous example corresponding to the decay model, as the locally structurally observable state $x$ corresponds exactly to the observed output of the linear model in Eq. \eqref{eq:output_linear}. Additionally, the product \(b z\) defines a \textit{locally structurally observable parameter-state combination} that does not correspond to the observed output. This finding highlights the fact that the combined local structural analysis based on parameter-state symmetries can find locally structurally observable parameter-state combinations that are not observed in the outputs.

To illustrate the parameter-state transformations by our class of parameter-state symmetries explicitly, we select \(\chi_{b}(\boldsymbol{\theta}) = 1\), which yields the generating vector field
\begin{equation}
    X = b \partial_{b} - z \partial_{z}\,.
    \label{eq:X_linear_particular}
\end{equation}
The corresponding parameter-state symmetry is then
\begin{equation}
    \Gamma^{\mathbf{x},\boldsymbol{\theta}}_{\varepsilon} : (t, x, z, a, b, c) \mapsto \left( t, x, z e^{-\varepsilon}, a, b e^{\varepsilon}, c \right)\,,
    \label{eq:parameter_symmetry_linear}
\end{equation}
which represents a scaling of the parameter \(b\) and the state \(z\) in opposite directions while preserving the model structure. Indeed, the model structure is maintained under this transformation, as can be verified by
\begin{equation}
    \begin{split}
        \widehat{\frac{\mathrm{d}x}{\mathrm{d}t}}(\varepsilon) &= \hat{a}(\varepsilon) \hat{x}(\varepsilon) + \hat{b}(\varepsilon) \hat{z}(\varepsilon)\\
        \widehat{\frac{\mathrm{d}z}{\mathrm{d}t}}(\varepsilon) &= \hat{c}(\varepsilon) \hat{z}(\varepsilon)\\
    \end{split}\quad \Longrightarrow \quad
    \begin{split}
      \frac{\mathrm{d}x}{\mathrm{d}t} &= a x + b e^{\varepsilon} z e^{-\varepsilon}\\
\frac{\mathrm{d}z}{\mathrm{d}t} e^{-\varepsilon} &= c z e^{-\varepsilon}      
     \end{split}\quad \Longrightarrow \quad
    \begin{split}
      \dfrac{\mathrm{d}x}{\mathrm{d}t}&= a x + b z\\
      \frac{\mathrm{d}z}{\mathrm{d}t} &= c z\\
    \end{split}\,.
    \label{eq:linear_sys_preserved}
\end{equation}
Thus, Eq.~\eqref{eq:linear_sys_preserved} explicitly demonstrates that the parameter-state symmetry satisfies the symmetry conditions in Eq.~\eqref{eq:sym_cond} in Section \ref{ssec:classical_symmetry}. This example, in line with the decoupled decay model, shows how parameter-state symmetries reveal locally structurally identifiable parameters and locally structurally observable states and parameter-state combinations.

\subsubsection{A non-autonomous glucose-insulin model}
We consider a model for the regulation of glucose and insulin in the blood originally studied by means of full symmetries in~\cite{massonis2020finding}. This model has two states, \(x_{1}(t)\) (glucose) and \(x_{2}(t)\) (insulin), described by

\begin{equation}
    \begin{split}
        \dfrac{\mathrm{d}x_{1}}{\mathrm{d}t}&=u+p_{1}x_{1}-p_{2}x_{2}\,,\\
        \dfrac{\mathrm{d}x_{2}}{\mathrm{d}t}&=p_{3}x_{2}+p_{4}x_{1}\,,\\        
    \end{split}
    \label{eq:linear_glu}
\end{equation}
where $u(t)$ is a time-dependent input that corresponds to the glucose that enters the digestive system. The observed output is a glucose measurement:
\begin{equation}
    y=\dfrac{x_{1}}{V_{p}}\,,
    \label{eq:output_glu}
\end{equation}
and the parameter vector is \(\boldsymbol{\theta} = (p_{1}, p_{2}, p_{3}, p_{4},V_{p})^\top\).

The family of generating vector fields for the glucose-insulin model is
\begin{equation}
    X=\chi_{p_{2}}(\boldsymbol{\theta})\left(p_{2}\partial_{p_{2}}-p_{4}\partial_{p_{4}}-x_{2}\partial_{x_{2}}\right)\,,
    \label{eq:X_glu_final}
\end{equation}
where $\chi_{p_{2}}:\mathbb{R}^{5}\to\mathbb{R}$ is an arbitrary function of the parameters. This vector field generates the full set of parameter-state symmetries for the glucose-insulin model. The universal invariants of the model are

\begin{equation}
    I_{1}=p_{1}\,,\quad{I}_{2}=p_{3}\,,\quad{I}_{3}=V_{p}\,,\quad{I}_{4}=p_{2}p_{4}\,,\quad{I}_{5}=x_{1}\,,\quad{I}_{6}=p_{2}x_{2}.
    \label{eq:uni_inv_glu}
\end{equation}
These universal invariants give a clear picture of the local structural properties of the glucose-insulin model. In terms of identifiability, the parameters $p_{1}$, $p_{3}$ and $V_{p}$ are locally structurally identifiable. Also, the parameters $p_{2}$ and $p_{4}$ are individually locally structurally \textit{unidentifiable}, whereas their product $p_{2}p_{4}$ is a locally structurally identifiable parameter combination.

In terms of observability, this example differs slightly from the previous decay model and linear model. For these previous models, the observed output constituted a structurally observable state or parameter-state combination. This is by definition always the case meaning that the observed output $x_{1}/V_{p}$ in Eq. \eqref{eq:output_glu} is also observable for the glucose-insulin model, but here we also have that the parameter $V_{p}$ is identifiable. This implies that the state $x_{1}$ in itself is locally structurally observable reflected by the universal invariant $I_{5}$ in Eq. \eqref{eq:uni_inv_glu}. Thus, in this case, unlike the previous examples, the locally structurally observable parameter-state combination given by the observed output is replaced by the locally structurally observable state $x_{1}$. In addition, we also have another structurally observable parameter-state combination given by the product $p_{2}x_{2}$. Accordingly, the glucose-insulin model also illustrates that locally structurally observable combinations may involve both states and parameters that are not included in observed outputs.

To illustrate the symmetry transformations explicitly, we select \(\chi_{p_{2}}(\boldsymbol{\theta})=1\), resulting in the vector field

\begin{equation}
    X=p_{2}\partial_{p_{2}}-p_{4}\partial_{p_{4}}-x_{2}\partial_{x_{2}}\,,
    \label{eq:X_glu_particular}
\end{equation}
which generates the parameter-state symmetry
\begin{equation}
    \Gamma^{\mathbf{x},\boldsymbol{\theta}}_{\varepsilon}:(t,x_{1},x_{2},p_{1},p_{2},p_{3},p_{4},V_{p})\mapsto\left(t,x_{1},x_{2}e^{-\varepsilon},p_{1},p_{2}e^{\varepsilon},p_{3},p_{4}e^{-\varepsilon},V_{p}\right)\,.
    \label{eq:parameter_symmetry_glu}
\end{equation}
This symmetry represents a scaling of the insulin state \(x_{2}\) and the parameters \(p_2\) and \(p_4\) in opposite directions, preserving the structure of the glucose-insulin model. Indeed, one can verify that

\begin{equation}
      \begin{split}
        \widehat{\dfrac{\mathrm{d}x_{1}}{\mathrm{d}t}}(\varepsilon)&=u+\hat{p}_{1}(\varepsilon)\hat{x}_{1}(\varepsilon)-\hat{p}_{2}(\varepsilon)\hat{x}_{2}(\varepsilon)\\
        \widehat{\dfrac{\mathrm{d}x_{2}}{\mathrm{d}t}}(\varepsilon)&=\hat{p}_{3}(\varepsilon)\hat{x}_{2}(\varepsilon)+\hat{p}_{4}(\varepsilon)\hat{x}_{1}(\varepsilon)\\        
    \end{split}\Longrightarrow
      \begin{split}
        \dfrac{\mathrm{d}x_{1}}{\mathrm{d}t}&=u+p_{1}x_{1}-p_{2}e^{\varepsilon}x_{2}e^{-\varepsilon}\\
        e^{-\varepsilon}\dfrac{\mathrm{d}x_{2}}{\mathrm{d}t}&=p_{3}x_{2}e^{-\varepsilon}+p_{4}e^{-\varepsilon}x_{1}\\        
    \end{split}\Longrightarrow
          \begin{split}
            \dfrac{\mathrm{d}x_{1}}{\mathrm{d}t}&=u+p_{1}x_{1}-p_{2}x_{2}\\
        \dfrac{\mathrm{d}x_{2}}{\mathrm{d}t}&=p_{3}x_{2}+p_{4}x_{1}\\        
    \end{split}\,,
    \label{eq:linear_glu_preserved}
\end{equation}
demonstrating that \(\Gamma^{\mathbf{x},\boldsymbol{\theta}}_{\varepsilon}\) satisfies the symmetry conditions in Eq.~\eqref{eq:sym_cond} of Section \ref{ssec:classical_symmetry}.

Notably, Massonis and Villaverde~\cite{massonis2020finding} independently identified this scaling-type symmetry using an ansatz-based polynomial approach for finding full symmetries. Thus, their work confirms that the parameter-state symmetries of the glucose-insulin model that we find.

\subsubsection{An epidemiological SEI model of the spread of tuberculosis}
We consider an epidemiological SEI model of the transmission of tuberculosis, originally analysed in~\cite{renardy2022structural}. To this end, let $S(t)$ be the density of susceptible individuals in a population, $E(t)$ be the corresponding exposed population and $I(t)$ be the population density for infected individuals. These states are described by
\begin{align}
  \dfrac{\mathrm{d}S}{\mathrm{d}t} &= - \beta I S + c - \mu_{S} S\,,\label{eq:ODE_S}\\
    \dfrac{\mathrm{d}E}{\mathrm{d}t} &= \beta \left(1 - \upsilon\right) I S - \delta E - \mu_{E} E\,,\label{eq:ODE_E}\\
    \dfrac{\mathrm{d}I}{\mathrm{d}t} &= \beta \upsilon I S + \delta E - \mu_{I} I\,.\label{eq:ODE_I}
\end{align}
Moreover, we observe the following two outputs
\begin{equation}
    y_{E}=k_{E}E\,,\quad{y}_{I}=k_{I}I\,,
    \label{eq:output_SEI}
\end{equation}
and their interpretation is that we observe the proportions $k_{E}$ and $k_{I}$ of the exposed and infected populations, respectively. In total, we have nine parameters collected in the vector $\boldsymbol{\theta}\in\mathbb{R}^{9}$ which are given by
\begin{equation}
\boldsymbol{\theta}=\begin{pmatrix}c\\\beta\\\mu_{S}\\\mu_{E}\\\mu_{I}\\\delta\\\upsilon\\k_{E}\\k_{I}\end{pmatrix}\,.
  \label{eq:param_SEI}
\end{equation}
The family of generating vector fields for the SEI model is given by:
\begin{equation}
  \begin{split}
    X &= -\left(\alpha_{1} + \alpha_{2} \upsilon\right) S\partial_{S} -\alpha_{1} E\partial_{E}-(\alpha_{1}+\alpha_{2}) I\partial_{I}\\
    &- c \left(\alpha_{1} + \alpha_{2} \upsilon\right)\partial_{c}+\beta \left(\alpha_{1} + \alpha_{2}\right)\partial_{\beta}+\alpha_{2} \delta\partial_{\mu_{E}}- \alpha_{2} \delta\partial_{\delta}\\
    &+\alpha_{2} \upsilon \left(\upsilon - 1\right)\partial_{\upsilon}+\alpha_{1}k_{E}\partial_{k_{E}}+(\alpha_{1}+\alpha_{2})k_{I}\partial_{k_{I}}\,,
\end{split}
  \label{eq:Lie_final}
\end{equation}
where $\alpha_{1},\alpha_{2}\in\mathbb{R}$ are two arbitrary parameters. These arbitrary parameters illustrate that we have two degrees of freedom in terms of constructing individual parameter-state symmetries for the SEI model compared to the corresponding one degree of freedom obtained in the previous three examples where we only could choose one arbitrary function. For this reason, we do not generate any particular parameter-state symmetry of the SEI model in order to demonstrate that it satisfies the symmetry conditions as for the previous examples.

The ten universal invariants of the parameter-state symmetries of the SEI model are given by
\begin{equation}
\begin{split}
    I_{1}&=\dfrac{S}{c}\,,\quad{I}_{2}=k_{E}E\,,\quad{I}_{3}=k_{I}I\,,\quad{I}_{4}=\mu_{S}\,,\quad{I}_{5}=\mu_{I}\,,\quad{I}_{6}=\mu_{E}+\delta\,,\\
    I_{7}&=\dfrac{k_{I}}{\beta}\,,\quad{I}_{8}=\delta\left(\dfrac{1-\upsilon}{\upsilon}\right)\,,\quad{I}_{9}=\beta{c}\upsilon\,,\quad{I}_{10}=\dfrac{\beta\delta}{k_{E}}\,.
\end{split}
    \label{eq:uni_inv_SEI}
\end{equation}
Giving their interpretation in reverse order, the seven universal invariants from $I_{4}$ to $I_{10}$ correspond to the locally structurally identifiable parameter combinations of the SEI model. Importantly, these universal parameter invariants of our parameter-state symmetries are the same as the universal parameter invariants of the parameter symmetries of the output system of the SEI model analysed in our previous work~\cite{borgqvist2024framingglobalstructuralidentifiability}. Moreover, these seven universal parameter invariants do not only correspond to the locally structurally identifiable parameter combinations, but in fact they are also globally structurally identifiable as they were the outcome of the standard differential algebra approach implemented in~\cite{renardy2022structural}. This shows that our current local structural identifiability analysis of the SEI model based on parameter-state symmetries is consistent with previous both local and global structural identifiability analyses based on distinct approaches. Better still, unlike those other approaches, our parameter-state symmetries also yield information about the structural observability of the SEI model. 

The universal invariants $I_{1}$, $I_{2}$ and $I_{3}$ in Eq. \eqref{eq:uni_inv_SEI} correspond to the locally structurally observable parameter-state combinations. In analogy with the decay model and the linear model, the latter two of these universal invariants are exactly the observed outputs of the SEI model in Eq.~\eqref{eq:output_SEI}. Moreover, the first universal invariant $I_{1}$ is a locally structurally observable parameter-state in terms of the unobserved susceptible population $S$. Again, the locally structurally observable parameter-state combination $S/c$ reveals that, similarly to the glucose-insulin model, observable parameter-state combinations can indeed be composed of parameters and states that are not included in the observed outputs.


\section{Discussion}
We have shown that the local structural identifiability and local structural observability of a model can be understood in terms of its universal invariants. This approach extends our previous work on parameter symmetries, which focused only on the output system depending solely on the observed outputs $\mathbf{y}$ and where the states $\mathbf{x}$ had been eliminated~\cite{borgqvist2024framingglobalstructuralidentifiability}. Here, we show that these parameter symmetries of the output system naturally generalise to parameter-state symmetries a full model as in Eq. \eqref{eq:ODE_sys_full}, which consists of a state-dependent ODE system together with equations for observed outputs. It is well known that this full model and the corresponding output system are output-equivalent~\cite{eisenberg2019inputoutputequivalenceidentifiabilitysimple}, meaning that observed outputs generated from one system also solve the other. Consequently, the parameter symmetries of the output system correspond directly to the parts of the parameter-state symmetries of the full model that depend on parameters, and all universal parameter invariants in both formulations are, in fact, \textit{identical}. Thus, just as the full model and its output system are output-equivalent, a corresponding equivalence exists between the parameter symmetries of the output system and the parameter-state symmetries of the full model. Better still, by analysing these parameter-state symmetries of the full model, we obtain additional information about local structural observability owing to the inclusion of states in the analysis. Accordingly, parameter-state symmetries afford a complete local structural analysis of any model, capturing both identifiability of parameters and observability of states.

There are several advantages to analysing the structural properties of models using the generalised parameter-state symmetries presented here compared with the output-dependent parameter symmetries introduced previously~\cite{borgqvist2024framingglobalstructuralidentifiability}. Most notably, by including the states in the analysis, parameter-state symmetries provide information not only about local structural identifiability of parameters but also about local structural observability of states. This means that the local structural observability of a model can be analysed by calculating its universal state invariants. Moreover, our approach eliminates the need to rewrite the model, composed of a state-dependent ODE system plus equations for observed outputs as in Eq.~\eqref{eq:ODE_sys_full}, as an equivalent output system. Within the standard differential algebra approach~\cite{ljung1994global}, this rewriting required computationally intensive algebraic manipulations and was restricted to underlying ODEs and observed outputs defined by rational functions of the states, inputs, and parameters. In contrast, the parameter-state symmetry analysis can, in principle, be applied to a wider class of models, without assumptions on the functional forms of the underlying ODEs or the observed outputs. While solving the linearised symmetry conditions may become algebraically complex for large or highly non-linear models, these difficulties are practical rather than theoretical, and numerical or computer-assisted approaches could be employed to obtain the parameter-state symmetries.

An interesting direction for future research is to explore the practical solvability of the linearised symmetry conditions for different types of state-dependent ODE systems and observed outputs. Technically, the parameter-state symmetries are obtained by solving $n$ linearised symmetry conditions for $n$ state-infinitesimals $\boldsymbol{\eta}$ and $p$ parameter-infinitesimals $\boldsymbol{\chi}$. Invariance of the observed outputs $\mathbf{y}$ imposes $m$ additional conditions, where $1<m<n$, on these infinitesimals. Some of the state-infinitesimals $\eta_{i}$, where $i\in\{1,\ldots,n\}$, and parameter-infinitesimals $\chi_{\ell}$, where $\ell\in\{1,\ldots,p\}$, can then be determined from these conditions. These expressions can subsequently be substituted back into the linearised symmetry conditions to solve for the remaining infinitesimals. While the analysis presented here demonstrates that the approach is generally applicable, there is no guarantee that the conditions can be solved in closed form for all models, particularly when the number of states, outputs, or parameters is large, or when the underlying ODEs are highly non-linear. Nevertheless, practical difficulties could be addressed using numerical or computer-assisted methods, including symbolic computation or machine learning-based approaches. A key question for future research is therefore: for which classes of state-dependent ODE systems and observed outputs is it feasible to determine parameter-state symmetries and thereby obtain information about local structural identifiability and local structural observability?

Building on the challenges outlined above, a natural next step is to develop automated approaches for performing combined local structural identifiability and local structural observability analyses based on parameter-state symmetries. For larger systems with many states, outputs, and parameters, manual calculations quickly become infeasible, making computer-assisted implementations essential. Encouragingly, inspiration can be drawn from existing automated frameworks for identifiability analyses, such as the Julia implementation of the differential algebra approach~\cite{dong2023structuralIdentifiabilityJL} and the Mathematica implementation of the EAR approach~\cite{anguelova2012minimal}, which successfully handle symbolic calculations for moderately complex systems. A desirable future direction is therefore to integrate parameter-state symmetries into such symbolic or numerical platforms, providing a practical tool for mechanistic modelling. Yet, a key open question remains: for which classes of state-dependent ODE systems, and with how many states, outputs, and parameters, can automated parameter-state-symmetry approaches provide efficient and reliable structural identifiability and structural observability analyses? Answering this question would constitute an important stepping-stone toward fully exploiting parameter-state symmetries for combined local structural identifiability and local structural observability analyses of mechanistic models in practice.




\section*{Data availability statement}
Regarding the combined SI- and SO- analysis of the SEI-model, details and relevant scripts are available at the public github-repository associated with this project; \url{https://github.com/JohannesBorgqvist/parameter_symmetries_identifiability_and_observability/}.

\section*{Acknowledgements}
JGB is funded by a grant from the Wenner-Gren foundations (Grant number: FT2023-0005). JBG thanks the Wenner-Gren foundations for a research fellowship. APB thanks the Mathematical Institute, Oxford for a Hooke Research Fellowship. This work was supported by a grant from the Simons Foundation (MP-SIP-00001828, REB)

\section*{CRediT author statment}
\begin{itemize}
\item[\textbf{JGB}] Conceptualization, Methodology, Visualization (made the figures), Writing - Original Draft, Writing - Review \& Editing, Formal analysis (derived and proved theorems and conducted calculations), Software (wrote scripts that conducted calculations). 
\item[\textbf{APB}] Conceptualization, Writing - Original Draft, Writing - Review \& Editing.
\item[\textbf{FO}] Conceptualization, Writing - Original Draft, Writing - Review \& Editing,
\item[\textbf{REB}] Conceptualization, Writing - Original Draft, Writing - Review \& Editing.
\end{itemize}

\clearpage
\appendix  

\renewcommand{\thesection}{\Alph{section}}

\counterwithin{equation}{section}  

\counterwithin{figure}{section}
\counterwithin{table}{section}

\section*{Appendices}


\section{A toy example of two decoupled decay equations}
We study the following decoupled system of decay ODEs
\begin{equation}
    \begin{split}
        \dfrac{\mathrm{d}u}{\mathrm{d}t}&=\kappa_{1}-\lambda{u}\,,\\
        \dfrac{\mathrm{d}v}{\mathrm{d}t}&=\kappa_{1}-\lambda{v}\,,\\        
    \end{split}
    \label{eq:decay_sys_app}
\end{equation}
together with the following observed output
\begin{equation}
    y=u+v\,.
    \label{eq:output_decay_app}
\end{equation}
In total, we have three parameters given by
\begin{equation}
    \boldsymbol{\theta}=\begin{pmatrix}\kappa_{1}\\\kappa_{2}\\\lambda\end{pmatrix}\,,
\end{equation}
where $\lambda$ is the common decay rate and where $\kappa_{1},\kappa_{2}$ are two constant formation rates. Given this system of equations, we want to find a family of infinitesimal generators of the Lie group of the form
\begin{equation}
        X=\eta_{u}(u,v,\boldsymbol{\theta})\partial_{u}+\eta_{v}(u,v,\boldsymbol{\theta})\partial_{v}
        +\chi_{\kappa_{1}}(\boldsymbol{\theta})\partial_{\kappa_{1}}+\chi_{\kappa_{2}}(\boldsymbol{\theta})\partial_{\kappa_{2}}+\chi_{\lambda}(\boldsymbol{\theta})\partial_{\lambda}\,.
    \label{eq:X_decay_general}
\end{equation}
The prolonged infinitesimals are given by
\begin{equation}
\begin{split}
    \eta_{u}^{(1)}&=D_{t}\eta_{u}=\dfrac{\mathrm{d}u}{\mathrm{d}t}\dfrac{\partial\eta_{u}}{\partial_{u}}+\dfrac{\mathrm{d}v}{\mathrm{d}t}\dfrac{\partial\eta_{u}}{\partial_{v}}\,,\\
    \eta_{v}^{(1)}&=D_{t}\eta_{v}=\dfrac{\mathrm{d}u}{\mathrm{d}t}\dfrac{\partial\eta_{v}}{\partial_{u}}+\dfrac{\mathrm{d}v}{\mathrm{d}t}\dfrac{\partial\eta_{v}}{\partial_{v}}\,,\\  
    \end{split}
    \label{eq:prolonged_decay}
\end{equation}
and the linearised symmetry conditions are given by
\begin{equation}
     \begin{split}
         \eta^{(1)}_{u}&=-\lambda\eta_{u}+\chi_{\kappa_{1}}-\chi_{\lambda}u\\
        \eta^{(1)}_{v}&=-\lambda\eta_{v}+\chi_{\kappa_{2}}-\chi_{\lambda}v         
     \end{split}\quad\quad\textrm{whenever}\quad\quad
  \begin{split}
        \dfrac{\mathrm{d}u}{\mathrm{d}t}&=\kappa_{1}-\lambda{u}\\
        \dfrac{\mathrm{d}v}{\mathrm{d}t}&=\kappa_{1}-\lambda{v}\\        
    \end{split}\,.    
     \label{eq:lin_sym_decoupled}
\end{equation}
Since the output $y$ in Eq. \eqref{eq:output_decay_app} must be invariant under transformations by these parameter symmetries, we impose the invariance condition $X(y)=0$ on the infinitesimal generators in Eq. \eqref{eq:X_decay_general}. This condition amounts to the following equation for the infinitesimals $\eta_{u}$ and $\eta_{v}$:
\begin{equation}
    \eta_{u}=-\eta_{v}\,.
    \label{eq:eta_eq_decay}
\end{equation}
Substituting this equation into the linearised symmetry conditions in Eq. \eqref{eq:lin_sym_decoupled} yields
\begin{align}
    \eta^{(1)}_{u}&=-\lambda\eta_{u}+\chi_{\kappa_{1}}-\chi_{\lambda}u\label{eq:eta_u_prolonged}\,,\\
    -\eta^{(1)}_{u}&=\lambda\eta_{u}+\chi_{\kappa_{2}}-\chi_{\lambda}v\,,
\end{align}
taking the sum of these two equations results in
\begin{equation}
    0=(\chi_{\kappa_{1}}+\chi_{\kappa_{2}})-\chi_{\lambda}(u+v)\,.
    \label{eq:decay_infinitesimals_eq}
\end{equation}
By decomposing this equation with respect to the linearly independent factors $\{1,(u+v)\}$, we obtain the following two equations
\begin{equation}
    \chi_{\kappa_{1}}=-\chi_{\kappa_{2}}\,,\quad\chi_{\lambda}=0\,,
    \label{eq:chi_decay_final}
\end{equation}
for the parameter infinitesimals. Substituting these parameter infinitesimals back into Eq. \eqref{eq:eta_u_prolonged}, we see that the infinitesimal $\eta_{u}$ solves the following PDE
\begin{equation}
    \chi_{\kappa_{1}}=D_{t}\eta_{u}-\lambda\eta_{u}\,. 
    \label{eq:eta_u_decay_final_app}
\end{equation}
Thus, the family of infinitesimal generators of the decay model is given by
\begin{equation}
\begin{split}
        X&=\eta_{u}(u,v,\boldsymbol{\theta})(\partial_{u}-\partial_{v})+\chi_{\kappa_{1}}(\boldsymbol{\theta})(\partial_{\kappa_{1}}-\partial_{\kappa_{2}})\,,
    \end{split}    
    \label{eq:X_decay_final_app}
\end{equation}
for some arbitrary function $\chi_{\kappa_{1}}:\mathbb{R}^{3}\mapsto\mathbb{R}$ and where $\eta_{u}$ solves Eq. \eqref{eq:eta_u_decay_final_app}. 

\subsection{Calculating universal invariants for the decay model}\label{ssec:invariant_decay}
A universal invariant of the decay model is a non-constant function $I=I(u(t),v(t),\boldsymbol{\theta})$ satisfying $X(I)=0$ where $X$ is given by Eq. \eqref{eq:X_decay_final_app}. Accordingly, $I$ solves the following linear PDE  
\begin{equation}
    0=\eta_{u}(u,v,\boldsymbol{\theta})\left(\dfrac{\partial I}{\partial_{u}}-\dfrac{\partial I}{\partial_{v}}\right)+\chi_{\kappa_{1}}(\boldsymbol{\theta})\left(\dfrac{\partial I}{\partial_{\kappa_{1}}}-\dfrac{\partial I}{\partial_{\kappa_{2}}}\right)\,,
    \label{eq:invariant_equation_decay}
\end{equation}
which we solve using the method of characteristics. If $s$ is a parameter that parametrises a solution curve $I(s)=I(u(s),v(s),\boldsymbol{\theta}(s))$, then the chain rule yields
\begin{equation}
    \dfrac{\mathrm{d}I}{\mathrm{d}s}=\dfrac{\mathrm{d}u}{\mathrm{d}s}\dfrac{\partial{I}}{\partial{u}}+\dfrac{\mathrm{d}v}{\mathrm{d}s}\dfrac{\partial{I}}{\partial{v}}+\dfrac{\mathrm{d}\kappa_{1}}{\mathrm{d}s}\dfrac{\partial{I}}{\partial{\kappa_{1}}}+\dfrac{\mathrm{d}\kappa_{2}}{\mathrm{d}s}\dfrac{\partial{I}}{\partial{\kappa_{2}}}+\dfrac{\mathrm{d}\lambda}{\mathrm{d}s}\dfrac{\partial{I}}{\partial{\lambda}}\,.  
    \label{eq:chain_rule_decay}
\end{equation}
By comparing Eqs. \eqref{eq:invariant_equation_decay} and \eqref{eq:chain_rule_decay}, we obtain the following characteristic equations:
\begin{align}
    \dfrac{\mathrm{d}I}{\mathrm{d}s}&=0\,,\label{eq:I_is_0}\\
    \dfrac{\mathrm{d}u}{\mathrm{d}s}&=\eta_{u}(u,v,\boldsymbol{\theta})\,,\label{eq:chara_u}\\
\dfrac{\mathrm{d}v}{\mathrm{d}s}&=-\eta_{u}(u,v,\boldsymbol{\theta})\,,\label{eq:chara_v}\\
    \dfrac{\mathrm{d}\kappa_{1}}{\mathrm{d}s}&=\chi_{\kappa_{1}}(\boldsymbol{\theta})\,,\label{eq:chara_k1}\\
    \dfrac{\mathrm{d}\kappa_{2}}{\mathrm{d}s}&=-\chi_{\kappa_{1}}(\boldsymbol{\theta})\,,\label{eq:chara_k2}\\  
\dfrac{\mathrm{d}\lambda}{\mathrm{d}s}&=0\,.\label{eq:chara_lambda}
\end{align}
The first equation in Eq. \eqref{eq:I_is_0} states that any invariant $I$ is an integration constant, i.e. $I=K$, or a so called first integral of the characteristic equations above. Consequently, Eq. \eqref{eq:chara_lambda} implies that $I_{1}=\lambda$ is a universal parameter invariant. By combining the characteristic equation of Eqs. \eqref{eq:chara_k1} and \eqref{eq:chara_k2}, we obtain the following
\begin{equation}
    \dfrac{\mathrm{d}\kappa_{1}}{\mathrm{d}\kappa_{2}}=-1\,,
    \label{eq:chara_k1_and_k2_decay}
\end{equation}
which is easily integrated and yields the second universal parameter invariant $I_{2}=\kappa_{1}+\kappa_{2}$. Similarly, the combination of the characteristic equations in Eqs. \eqref{eq:chara_u} and \eqref{eq:chara_v} results in
\begin{equation}
    \dfrac{\mathrm{d}u}{\mathrm{d}v}=-1\,,
    \label{eq:chara_u_and_v_decay}
\end{equation}
which can be integrated and gives the universal invariant $I_{3}=u+v$. In summary, the three universal invariants of the decay model are given by
\begin{equation}
    I_{1}=\lambda\,,\quad{I}_{2}=\kappa_{1}+\kappa_{2}\,,\quad{I}_{3}=u+v\,.
    \label{eq:uni_inv_decay_app}
\end{equation}
\subsection{Generating a parameter symmetry for the decay model}
Next, we study the particular parameter symmetry of the decay model of interest defined by choosing the infinitesimal for $u$ to a constant. Under this assumption, it follows that $D_{t}\eta_{u}=0$, and by Eq. \eqref{eq:eta_u_decay_final_app} such a constant infinitesimal is given by
\begin{equation}
    \eta_{u}=-\dfrac{\chi_{\kappa_{1}}}{\lambda}\,.
    \label{eq:eta_u_decay_particular}
\end{equation}
For the purpose of illustrating a specific parameter symmetry of the decay model, we also set the value of the parameter infinitesimal $\chi_{\kappa_{1}}$ to one, i.e. $\chi_{\kappa_{1}}=1$, which gives the particular generating vector field
\begin{equation}
\begin{split}
        X&=\dfrac{1}{\lambda}(\partial_{u}-\partial_{v})+(\partial_{\kappa_{1}}-\partial_{\kappa_{2}})\,.
    \end{split}    
    \label{eq:X_decay_particular_app}
\end{equation}
The parameter symmetry $\Gamma_{\varepsilon}^{\mathbf{x},\boldsymbol{\theta}}:(t,u,v,\kappa_{1},\kappa_{2},\lambda)\mapsto(t,\hat{u}(\varepsilon),\hat{v}(\varepsilon),\hat{\kappa_{1}}(\varepsilon),\hat{\kappa_{2}}(\varepsilon),\hat{\lambda}(\varepsilon))$ that is generated by $X$ in Eq. \eqref{eq:X_decay_particular_app} solves the following system of ODEs
\begin{align}
    \dfrac{\mathrm{d}\hat{u}}{\mathrm{d}\varepsilon}&=\dfrac{1}{\lambda}\,,\quad\hat{u}(\varepsilon=0)=u\,,\label{eq:u_trans_decay}\\
\dfrac{\mathrm{d}\hat{v}}{\mathrm{d}\varepsilon}&=-\dfrac{1}{\lambda}\,,\quad\hat{v}(\varepsilon=0)=v\,,\label{eq:v_trans_decay}\\ 
    \dfrac{\mathrm{d}\hat{\kappa_{1}}}{\mathrm{d}\varepsilon}&=1\,,\quad\hat{\kappa_{1}}(\varepsilon=0)=\kappa_{1}\,,\label{eq:k1_trans_decay}\\
\dfrac{\mathrm{d}\hat{\kappa_{2}}}{\mathrm{d}\varepsilon}&=-1\,,\quad\hat{\kappa_{2}}(\varepsilon=0)=\kappa_{2}\,,
\label{eq:k2_trans_decay}\\    
\dfrac{\mathrm{d}\hat{\lambda}}{\mathrm{d}\varepsilon}&=0\,,\quad\hat{\lambda}(\varepsilon=0)=\lambda\,.\label{eq:lambda_trans_decay}
\end{align}
By solving these ODEs, the resulting parameter symmetry is given by
\begin{equation}
    \Gamma_{\varepsilon}^{\mathbf{x},\boldsymbol{\theta}}:(t,u,v,\kappa_{1},\kappa_{2},\lambda)\mapsto\left(t,u+\dfrac{\varepsilon}{\lambda},v-\dfrac{\varepsilon}{\lambda},\kappa_{1}+\varepsilon,\kappa_{2}-\varepsilon,\lambda\right)\,.
    \label{eq:parameter_symmetry_sum_app}
\end{equation}
As expected, transformations by the parameter symmetry $\Gamma_{\varepsilon}^{\mathbf{x},\boldsymbol{\theta}}$ in Eq. \eqref{eq:parameter_symmetry_sum_app} preserve the model structure. This is illustrated by the following calculations:
\begin{equation}
    \begin{split}
        \widehat{\dfrac{\mathrm{d}u}{\mathrm{d}t}}(\varepsilon)&=\hat{\kappa_{1}}(\varepsilon)-\hat{\lambda}(\varepsilon){\hat{u}(\varepsilon)}\\
        \widehat{\dfrac{\mathrm{d}v}{\mathrm{d}t}}(\varepsilon)&=\hat{\kappa_{2}}(\varepsilon)-\hat{\lambda}(\varepsilon){\hat{v}(\varepsilon)}\\        
    \end{split}\quad\Longrightarrow\quad\begin{split}
        \dfrac{\mathrm{d}u}{\mathrm{d}t}&=(\kappa_{1}+\varepsilon)-\lambda{\left(u+\dfrac{\varepsilon}{\lambda}\right)}\\
        \dfrac{\mathrm{d}v}{\mathrm{d}t}&=(\kappa_{2}-\varepsilon)-\lambda{\left(v-\dfrac{\varepsilon}{\lambda}\right)}\\        
    \end{split}\quad\Longrightarrow\quad
        \begin{split}
        \dfrac{\mathrm{d}u}{\mathrm{d}t}&=\kappa_{1}-\lambda{u}\\
        \dfrac{\mathrm{d}v}{\mathrm{d}t}&=\kappa_{1}-\lambda{v}\\        
    \end{split}\,.
    \label{eq:decay_sys_preserved_app}
\end{equation}

\section{A linear model}
Next, we study the following linear system
\begin{equation}
    \begin{split}
        \dfrac{\mathrm{d}x}{\mathrm{d}t}&=ax+bz\,,\\
        \dfrac{\mathrm{d}z}{\mathrm{d}t}&=cz\,,\\        
    \end{split}
    \label{eq:linear_sys_app}
\end{equation}
together with the following observed output
\begin{equation}
    y=x\,,
    \label{eq:output_linear_app}
\end{equation}
and the following parameters;
\begin{equation}
    \boldsymbol{\theta}=\begin{pmatrix}a\\b\\c\end{pmatrix}\,.
    \label{eq:param_linear}
\end{equation}
We are looking for a family of infinitesimal generators of the Lie group of the form
\begin{equation}
    \begin{split}
        X&=\eta_{u}(x,z,\boldsymbol{\theta})\partial_{x}+\eta_{z}(x,z,\boldsymbol{\theta})\partial_{z}+\chi_{a}(\boldsymbol{\theta})\partial_{a}+\chi_{b}(\boldsymbol{\theta})\partial_{b}+\chi_{c}(\boldsymbol{\theta})\partial_{c}\,.
    \end{split}
    \label{eq:X_linear}
\end{equation}
The linearised symmetry conditions are given by
\begin{equation}
     \begin{split}
         \eta^{(1)}_{x}&=a\eta_{x}+\chi_{a}x+b\eta_{z}+\chi_{b}z\\
        \eta^{(1)}_{z}&=c\eta_{z}+\chi_{c}z        
     \end{split}\quad\quad\textrm{whenever}\quad\quad
  \begin{split}
        \dfrac{\mathrm{d}x}{\mathrm{d}t}&=ax+bz\\
        \dfrac{\mathrm{d}z}{\mathrm{d}t}&=cz\\        
    \end{split}\,.    
     \label{eq:lin_sym_linear}
\end{equation}
The output $y$ in Eq. \eqref{eq:output_linear_app} is invariant according to the equation $X(y)=0$ which implies that
\begin{equation}
     \eta_{x}=0\,,
    \label{eq:eta_x_linear}
\end{equation}
and thus the infinitesimal for $x$ is zero. Accordingly, we also have $D_{t}\eta_{x}$ which means that from the linearised symmetry condition for $x$ in Eq. \eqref{eq:lin_sym_linear}, we obtain
\begin{equation}
    \eta_{z}=-\dfrac{1}{b}\left(\chi_{a}x+\chi_{b}z\right)\,.
    \label{eq:eta_z_temp_1}
\end{equation}
By substituting Eq. \eqref{eq:eta_z_temp_1} into the second linearised symmetry condition for $z$ in Eq. \eqref{eq:lin_sym_linear}, we obtain
\begin{equation}
    -\dfrac{\chi_{a}}{b}(ax+bz)-\dfrac{\chi_{b}}{b}(cz)=-\dfrac{c}{b}\left(\chi_{a}x+\chi_{b}z\right)+\chi_{c}z\,,
    \label{eq:z_lin_sym_1}
\end{equation}
which can be simplified to
\begin{equation}
    0=\chi_{a}\left(\dfrac{c-a}{b}\right)x+(\chi_{a}-\chi_{c})z\,.
    \label{eq:z_lin_sym_2}
\end{equation}
By decomposing the above equation with respect to the linearly independent set $\{x,z\}$, we obtain that the only two solutions to the above equation are $\chi_{a}=\chi_{c}=0$. Substituting these solutions back into Eq. \eqref{eq:eta_z_temp}, we obtain
\begin{equation}
    \eta_{z}=-\left(\dfrac{\chi_{b}}{b}\right)z\,.
    \label{eq:eta_z_temp}
\end{equation}
Hence, the family of generating vector fields for the linear model can be written in the form\footnote{The final vector field $X$ in Eq. \eqref{eq:X_linear_final_app} for the linear model was obtained in three steps: First, we substitute the expression for $\eta_{z}$ in Eq. \eqref{eq:eta_z_temp} into the general equation for $X$ in Eq. \eqref{eq:X_linear} together with the substitutions $\chi_{a}=\chi_{c}=\eta_{x}=0$. Second, we factor out $\chi_{b}/b$. Third, we rename the arbitrary function according to $\chi_{b}/b\mapsto\chi_{b}$.}
\begin{equation}
    X=\chi_{b}(\boldsymbol{\theta})(b\partial_{b}-z\partial_{z})\,,
    \label{eq:X_linear_final_app}
\end{equation}
for an arbitrary function $\chi_{b}:\mathbb{R}^{3}\mapsto\mathbb{R}$ of the parameters $\boldsymbol{\theta}$ in Eq. \eqref{eq:param_linear}.

\subsection{Calculating universal invariants for the linear model}
We have already seen that $\chi_{a}=\chi_{c}=0$. Hence, the parameters $I_{1}=a$ and $I_{2}=c$ are universal parameter invariants and therefore structurally locally identifiable. Also, $\eta_{x}=0$ and hence the state $I_{3}=x$ is locally structurally observable, which is hardly surprising given that it corresponds to the observed output $y$ in Eq. \eqref{eq:output_linear_app}. For the remaining universal invariant, the method of characteristics applied on the generator $X$ in Eq. \eqref{eq:X_linear_final_app} yields the following characteristic equation
\begin{equation}
    \dfrac{\mathrm{d}z}{\mathrm{d}b}=-\dfrac{z}{b}\,,
    \label{eq:linear_uni_inv_final }
\end{equation}
which can be integrated to the following first integral $I_{4}=bz$. In summary, the four universal invariants of the linear model are given by
\begin{equation}
    I_{1}=a\,,\quad{I}_{2}=c\,,\quad{I}_{3}=x\,,\quad{I}_{4}=bz.
    \label{eq:uni_inv_linear_app}
\end{equation}
\subsection{Generating a parameter symmetry for the linear model}
For the choice $\chi_{b}(\theta)=1$, the corresponding generating vector field is given by 
\begin{equation}
    X=b\partial_{b}-z\partial_{z}\,.
    \label{eq:X_linear_particular_app}
\end{equation}
Generating the corresponding parameter symmetry, by using the same type of calculations as for the previously considered decay model, yields
\begin{equation}
    \Gamma_{\varepsilon}^{\mathbf{x},\boldsymbol{\theta}}:(t,x,z,a,b,c)\mapsto\left(t,x,ze^{-\varepsilon},a,be^{\varepsilon},c\right)\,.
    \label{eq:parameter_symmetry_linear_app}
\end{equation}
Again, transformations by the parameter symmetry $\Gamma_{\varepsilon}^{\mathbf{x},\boldsymbol{\theta}}$ in Eq. \eqref{eq:parameter_symmetry_linear_app} preserves the model structure of the linear model in Eq. \eqref{eq:linear_sys_app}. This is illustrated by the following calculations:
\begin{equation}
     \begin{split}
        \widehat{\dfrac{\mathrm{d}x}{\mathrm{d}t}}(\varepsilon)&=\hat{a}(\varepsilon)\hat{x}(\varepsilon)+\hat{b}(\varepsilon)\hat{z}(\varepsilon)\\
        \widehat{\dfrac{\mathrm{d}z}{\mathrm{d}t}}(\varepsilon)&=\hat{c}(\varepsilon)\hat{z}(\varepsilon)\\        
    \end{split}\quad\Longrightarrow\quad
     \begin{split}
        \dfrac{\mathrm{d}x}{\mathrm{d}t}&=ax+be^{\varepsilon}ze^{-\varepsilon}\\
        \dfrac{\mathrm{d}z}{\mathrm{d}t}e^{-\varepsilon}&=cze^{-\varepsilon}\\        
    \end{split}\quad\Longrightarrow\quad
         \begin{split}
        \dfrac{\mathrm{d}x}{\mathrm{d}t}&=ax+bz\\
        \dfrac{\mathrm{d}z}{\mathrm{d}t}&=cz\\        
    \end{split}\,.
    \label{eq:linear_sys_preserved_app}
\end{equation}

\section{A non-autonomous glucose insulin model}
We consider a model for the regulation of glucose and insulin in the blood. Let $x_{1}(t)$ describe the glucose concentration in the blood at time $t$ and $x_{2}(t)$ describe the insulin concentration in the blood at time $t$. Then, the model is given by
\begin{equation}
    \begin{split}
        \dfrac{\mathrm{d}x_{1}}{\mathrm{d}t}&=u+p_{1}x_{1}-p_{2}x_{2}\,,\\
        \dfrac{\mathrm{d}x_{2}}{\mathrm{d}t}&=p_{3}x_{2}+p_{4}x_{1}\,,\\
    \end{split}
    \label{eq:linear_glu_app}
\end{equation}
where $u(t)$ is a time-dependent input that corresponds to the glucose that enters the digestive system. Moreover, the output corresponds to a glucose measurement.
\begin{equation}
    y=\dfrac{x_{1}}{V_{p}}\,,
    \label{eq:output_glu_app}
\end{equation}
and the model has five parameters:
\begin{equation}
    \boldsymbol{\theta}=\begin{pmatrix}p_{1}\\p_{2}\\p_{3}\\p_{4}\\V_{p}\end{pmatrix}\,.
    \label{eq:param_glu}
\end{equation}
Next, we are looking for a generating vector field $X$ of the following form
\begin{equation}
    \begin{split}
        X&=\eta_{1}(t,x_{1},x_{2},\boldsymbol{\theta})\partial_{x_{1}}+\eta_{2}(t,x_{1},x_{2},\boldsymbol{\theta})\partial_{x_{2}}+\chi_{p_{1}}(\boldsymbol{\theta})\partial_{p_{1}}+\chi_{p_{2}}(\boldsymbol{\theta})\partial_{p_{2}}\\
        &\quad+\chi_{p_{3}}(\boldsymbol{\theta})\partial_{p_{3}}+\chi_{p_{4}}(\boldsymbol{\theta})\partial_{p_{4}}+\chi_{V_{p}}(\boldsymbol{\theta})\partial_{V_{p}}\,.
    \end{split}
    \label{eq:X_glu_general}
\end{equation}
The linearised symmetry conditions are given by
\begin{equation}
    \begin{split}
        \dfrac{\mathrm{d}x_{1}}{\mathrm{d}t}\dfrac{\partial\eta_{1}}{\partial x_{1}}+\dfrac{\mathrm{d}x_{2}}{\mathrm{d}t}\dfrac{\partial\eta_{1}}{\partial x_{2}}&=p_{1}\eta_{1}-p_{2}\eta_{2}+\chi_{p_{1}}x_{1}-\chi_{p_{2}}x_{2}\\
        \dfrac{\mathrm{d}x_{1}}{\mathrm{d}t}\dfrac{\partial\eta_{2}}{\partial x_{1}}+\dfrac{\mathrm{d}x_{2}}{\mathrm{d}t}\dfrac{\partial\eta_{2}}{\partial x_{2}}&=p_{3}\eta_{2}+p_{4}\eta_{1}+\chi_{p_{3}}x_{2}+\chi_{p_{4}}x_{1}\\       \end{split}\quad\quad\textrm{whenever}\quad\quad
      \begin{split}
        \dfrac{\mathrm{d}x_{1}}{\mathrm{d}t}&=u+p_{1}x_{1}-p_{2}x_{2}\\
        \dfrac{\mathrm{d}x_{2}}{\mathrm{d}t}&=p_{3}x_{2}+p_{4}x_{1}\\        
    \end{split}\,.    
     \label{eq:lin_sym_glu}
\end{equation}
Since the output $y$ in Eq. \eqref{eq:output_glu_app} is invariant implying that $X(y)=0$ for $X$ in Eq. \eqref{eq:X_glu_general}, we obtain the following equation for the infinitesimal for $x_{1}$:
\begin{equation}
    \eta_{1}=\left(\dfrac{\chi_{V_{p}}}{V_{p}}\right)x_{1}\,.
    \label{eq:eta_1_temp}
\end{equation}
Substituting this equation into the first linearised symmetry condition for $x_{1}$ in Eq. \eqref{eq:lin_sym_glu}, yields the following equation
\begin{equation}
(u+p_{1}x_{1}-p_{2}x_{2})\left(\dfrac{\chi_{V_{p}}}{V_{p}}\right)=p_{1}\left(\dfrac{\chi_{V_{p}}}{V_{p}}\right)x_{1}-p_{2}\eta_{2}+\chi_{p_{1}}x_{1}-\chi_{p_{2}}x_{2}\,.
    \label{eq:lin_sym_x1}
\end{equation}
We note that the only term in the above equation that depends explicitly on time is the known input $u(t)$ on the left hand side. Since there is no time-dependence on the right hand side of Eq. \eqref{eq:lin_sym_x1}, the coefficient of $u(t)$ on the left hand side must be zero. From this follows that $\chi_{V_{p}}=0$ which, in turn, implies that $\eta_{1}=0$ according to Eq. \eqref{eq:eta_1_temp}. By substituting $\eta_{1}=0$ into Eq. \eqref{eq:lin_sym_x1} and solving the resulting equation for $\eta_{2}$ yields
\begin{equation}
    \eta_{2}=\dfrac{1}{p_{2}}\left(\chi_{p_{1}}x_{1}-\chi_{p_{2}}x_{2}\right)\,.
     \label{eq:eta_2_temp}
\end{equation}
Substituting $\eta_{2}$ in the second linearised symmetry condition for $x_{2}$ in Eq. \eqref{eq:lin_sym_glu} with the right hand side of Eq. \eqref{eq:eta_2_temp} yields
\begin{equation}
(u+p_{1}x_{1}-p_{2}x_{2})\left(\dfrac{\chi_{p_{1}}}{p_{2}}\right)-(p_{3}x_{2}+p_{4}x_{1})\left(\dfrac{\chi_{p_{2}}}{p_{2}}\right)=\dfrac{p_{3}}{p_{2}}\left(\chi_{p_{1}}x_{1}-\chi_{p_{2}}x_{2}\right)+\chi_{p_{3}}x_{2}+\chi_{p_{4}}x_{1}\,.
    \label{eq:lin_sym_x2_1}
\end{equation}
Using the same reasoning as in the case of the previous linearised symmetry condition, the coefficient of the time dependent input $u(t)$ is zero, and hence $\chi_{p_{1}}=0$. Substituting this value into Eq. \eqref{eq:lin_sym_x2_1} yields
\begin{equation}
0=p_{4}\left(\dfrac{\chi_{p_{2}}}{p_{2}}+\dfrac{\chi_{p_{4}}}{p_{4}}\right)x_{1}+\chi_{p_{3}}x_{2}\,.
    \label{eq:lin_sym_x2_2}
\end{equation}
Lastly, splitting up this equation with respect to the linearly independent set $\{x_{1},x_{2}\}$ we obtain that $\chi_{p_{3}}=0$ and $\chi_{p_{4}}=-(p_{4}/p_{2})\chi_{p_{2}}$. From Eq. \eqref{eq:eta_2_temp}, we have
\begin{equation}
    \eta_{2}=-\left(\dfrac{\chi_{p_{2}}}{p_{2}}\right)x_{2}\,.
     \label{eq:eta_2_final}
\end{equation}
The conclusion from these calculations is that the family of generating vector fields for the glucose-insulin model is given by
\begin{equation}
    X=\chi_{p_{2}}(\boldsymbol{\theta})\left(p_{2}\partial_{p_{2}}-p_{4}\partial_{p_{4}}-x_{2}\partial_{x_{2}}\right)\,,
    \label{eq:X_glu_final_app}
\end{equation}
where $\chi_{p_{2}}:\mathbb{R}^{5}\mapsto\mathbb{R}$ is an arbitrary function of the parameters $\boldsymbol{\theta}$ in Eq. \eqref{eq:param_glu}. 

\subsection{Calculating universal invariants for the glucose-insulin model}
We have already seen that $\chi_{p_{1}}=\chi_{p_{3}}=\chi_{V_{p}}=0$, and hence $I_{1}=p_{1}$, $I_{2}=p_{3}$ and $I_{3}=V_{p}$ are universal parameter invariants meaning that they are locally structurally identifiable. Similarly, $\eta_{1}=0$ implies that $I_{5}=x_{1}$ is a universal invariant and therefore $x_{1}$ is locally structurally observable. For the remaining two universal invariants, we apply the method of characteristics on the vector fields $X$ in Eq. \eqref{eq:X_glu_final_app}. Combining the characteristic equations for $p_{2}$ and $p_{4}$, we get
\begin{equation}
    \dfrac{\mathrm{d}p_{2}}{\mathrm{d}p_{4}}=-\dfrac{p_{2}}{p_{4}}\,,
    \label{eq:chara_p2_p4}
\end{equation}
which can be integrated and this results in the universal parameter invariant $I_{4}=p_{2}p_{4}$. In the same fashion, combining the characteristic equations for $x_{2}$ and $p_{2}$, yields
\begin{equation}
    \dfrac{\mathrm{d}x_{2}}{\mathrm{d}p_{2}}=-\dfrac{x_{2}}{p_{2}}\,,
    \label{eq:chara_x2_x2}
\end{equation}
which can be integrated, resulting in the universal invariant $I_{6}=p_{2}x_{2}$. In summary, the six universal invariants of the glucose-insulin model are given by
\begin{equation}
    I_{1}=p_{1}\,,\quad{I}_{2}=p_{3}\,,\quad{I}_{3}=V_{p}\,,\quad{I}_{4}=p_{2}p_{4}\,,\quad{I}_{5}=x_{1}\,,\quad{I}_{6}=p_{2}x_{2}.
    \label{eq:uni_inv_glu_app}
\end{equation}

\subsection{Generating a parameter symmetry for the glucose-insulin model}
The choice $\chi_{p_{2}}(\boldsymbol{\theta})=1$ results in the vector field
\begin{equation}
    X=p_{2}\partial_{p_{2}}-p_{4}\partial_{p_{4}}-x_{2}\partial_{x_{2}}\,,
    \label{eq:X_glu_particular_app}
\end{equation}
which generates the following parameter symmetry
\begin{equation}
    \Gamma_{\varepsilon}^{\mathbf{x},\boldsymbol{\theta}}:(t,x_{1},x_{2},p_{1},p_{2},p_{3},p_{4},V_{p})\mapsto\left(t,x_{1},x_{2}e^{-\varepsilon},p_{1},p_{2}e^{\varepsilon},p_{3},p_{4}e^{-\varepsilon},V_{p}\right)\,.
    \label{eq:parameter_symmetry_glu_app}
\end{equation}
As in the previous cases, transformations by the parameter symmetry $\Gamma_{\varepsilon}^{\mathbf{x},\boldsymbol{\theta}}$ in Eq. \eqref{eq:parameter_symmetry_glu_app} preserves the model structure of the glucose-insulin model in Eq. \eqref{eq:linear_glu_app}. This is illustrated by the following calculations:
\begin{equation}
      \begin{split}
        \widehat{\dfrac{\mathrm{d}x_{1}}{\mathrm{d}t}}(\varepsilon)&=u+\hat{p}_{1}(\varepsilon)\hat{x}_{1}(\varepsilon)-\hat{p}_{2}(\varepsilon)\hat{x}_{2}(\varepsilon)\\
        \widehat{\dfrac{\mathrm{d}x_{2}}{\mathrm{d}t}}(\varepsilon)&=\hat{p}_{3}(\varepsilon)\hat{x}_{2}(\varepsilon)+\hat{p}_{4}(\varepsilon)\hat{x}_{1}(\varepsilon)\\        
    \end{split}\Longrightarrow
      \begin{split}
        \dfrac{\mathrm{d}x_{1}}{\mathrm{d}t}&=u+p_{1}x_{1}-p_{2}e^{\varepsilon}x_{2}e^{-\varepsilon}\\
        e^{-\varepsilon}\dfrac{\mathrm{d}x_{2}}{\mathrm{d}t}&=p_{3}x_{2}e^{-\varepsilon}+p_{4}e^{-\varepsilon}x_{1}\\        
    \end{split}\Longrightarrow
          \begin{split}
        \dfrac{\mathrm{d}x_{1}}{\mathrm{d}t}&=u+p_{1}x_{1}-p_{2}x_{2}\\
        \dfrac{\mathrm{d}x_{2}}{\mathrm{d}t}&=p_{3}x_{2}+p_{4}x_{1}\\        
    \end{split}\,.
    \label{eq:linear_glu_preserved_app}
\end{equation}
\section{An epidemiological SEI model of the spread of tuberculosis}

We consider the following SEI model of epidemiological transmission of tuberculosis:
\begin{align}
  \dfrac{\mathrm{d}S}{\mathrm{d}t} &= - \beta I S + c - \mu_{S} S\,,\label{eq:ODE_S_app}\\
  \dfrac{\mathrm{d}E}{\mathrm{d}t} &= \beta \left(1 - \upsilon\right) I S - \delta E - \mu_{E} E\,,\label{eq:ODE_E_app}\\
  \dfrac{\mathrm{d}I}{\mathrm{d}t} &= \beta \upsilon I S + \delta E - \mu_{I} I\,.\label{eq:ODE_I_app}
\end{align}
We also observe the following two outputs
\begin{align}
  y_{E}&=k_{E}E\,,\label{eq:output_E_app}\\
  y_{I}&=k_{I}I\,,\label{eq:output_I_app}  
\end{align}
and their interpretation is that we observe the proportions $k_{E}$ and $k_{I}$ of the exposed and infected populations, respectively. In total, we have nine parameters collected in the vector $\boldsymbol{\theta}\in\mathbb{R}^{9}$ which are given by
\begin{equation}
\boldsymbol{\theta}=\begin{pmatrix}c\\\beta\\\mu_{S}\\\mu_{E}\\\mu_{I}\\\delta\\\upsilon\\k_{E}\\k_{I}\end{pmatrix}\,.
  \label{eq:param_SEI_app}
\end{equation}
We are looking for a family of infinitesimal generators of the Lie group given by
\begin{equation}
  \begin{split}
    X&=\eta_{S}(S,E,I,\boldsymbol{\theta})\partial_{S}+\eta_{E}(S,E,I,\boldsymbol{\theta})\partial_{E}+\eta_{I}(S,E,I,\boldsymbol{\theta})\partial_{I}\\
    &\quad+\chi_{c}(\boldsymbol{\theta})\partial_{c}+\chi_{\beta}(\boldsymbol{\theta})\partial_{\beta}+\chi_{\mu_{S}}(\boldsymbol{\theta})\partial_{\mu_{S}}+\chi_{\mu_{E}}(\boldsymbol{\theta})\partial_{\mu_{E}}+\chi_{\mu_{I}}(\boldsymbol{\theta})\partial_{\mu_{I}}+\chi_{\delta}(\boldsymbol{\theta})\partial_{\mu_{\delta}}\\
    &\quad+\chi_{\upsilon}(\boldsymbol{\theta})\partial_{\mu_{\upsilon}}+\chi_{k_{E}}(\boldsymbol{\theta})\partial_{k_{E}}+\chi_{k_{I}}(\boldsymbol{\theta})\partial_{k_{I}}\,.
    \end{split}
  \label{eq:X_SEI_original}
\end{equation}
The first three prolongations of the infinitesimals for the states are given by
\begin{align}
  \eta_{S}^{(1)}\left(S,E,I,\dfrac{\mathrm{d}S}{\mathrm{d}t},\boldsymbol{\theta}\right)&=D_{t}\eta_{S}(S,E,I,\boldsymbol{\theta})=\dfrac{\mathrm{d}S}{\mathrm{d}t}\dfrac{\partial\eta_{S}}{\partial S}+\dfrac{\mathrm{d}E}{\mathrm{d}t}\dfrac{\partial\eta_{S}}{\partial E}+\dfrac{\mathrm{d}I}{\mathrm{d}t}\dfrac{\partial\eta_{S}}{\partial I}\label{eq:eta_S_1}\,,\\
  \eta_{E}^{(1)}\left(S,E,I,\dfrac{\mathrm{d}E}{\mathrm{d}t},\boldsymbol{\theta}\right)&=D_{t}\eta_{E}(S,E,I,\boldsymbol{\theta})=\dfrac{\mathrm{d}S}{\mathrm{d}t}\dfrac{\partial\eta_{E}}{\partial S}+\dfrac{\mathrm{d}E}{\mathrm{d}t}\dfrac{\partial\eta_{E}}{\partial E}+\dfrac{\mathrm{d}I}{\mathrm{d}t}\dfrac{\partial\eta_{E}}{\partial I}\label{eq:eta_E_1}\,,\\
\eta_{I}^{(1)}\left(S,E,I,\dfrac{\mathrm{d}S}{\mathrm{d}t},\boldsymbol{\theta}\right)&=D_{t}\eta_{I}(S,E,I,\boldsymbol{\theta})=\dfrac{\mathrm{d}S}{\mathrm{d}t}\dfrac{\partial\eta_{I}}{\partial S}+\dfrac{\mathrm{d}E}{\mathrm{d}t}\dfrac{\partial\eta_{I}}{\partial E}+\dfrac{\mathrm{d}I}{\mathrm{d}t}\dfrac{\partial\eta_{I}}{\partial I}\label{eq:eta_I_1}\,,
\end{align}
where the total derivative is defined by: $D_{t}=\partial_{t}+(\mathrm{d}S/\mathrm{d}t)\partial_{S}+(\mathrm{d}E/\mathrm{d}t)\partial_{E}+(\mathrm{d}I/\mathrm{d}t)\partial_{I}$. These prolongations define the first prolongation of the infinitesimal generator $X^{(1)}$ according to
\begin{equation}
  X^{(1)}=X+\eta_{S}^{(1)}\partial_{\mathrm{d}S/\mathrm{d}t}+\eta_{E}^{(1)}\partial_{\mathrm{d}E/\mathrm{d}t}+\eta_{I}^{(1)}\partial_{\mathrm{d}I/\mathrm{d}t}\,.
  \label{eq:X_1_SEI}
\end{equation}
Before, we define the linearised symmetry conditions, we use the output invariance conditions to get conditions on two of the state infinitesimals. Specifically, the fact that the observed outputs are differential invariants of our generator yields equations for the infinitesimals $\eta_{E}$ and $\eta_{I}$, respectively. Starting with the infinitesimal for $E$, we have
$$X(y_{E})=0\Longrightarrow\quad{k_{E}}\eta_{E}+\chi_{K_{E}}E=0\,,$$
which gives us the following equation for $\eta_{E}$
\begin{equation}
  \eta_{E}=-\left(\dfrac{\chi_{k_{E}}}{k_{E}}\right)E\,.
  \label{eq:eta_E}
\end{equation}
Analogously, the equation for $\eta_{I}$ is given by
\begin{equation}
  \eta_{I}=-\left(\dfrac{\chi_{k_{I}}}{k_{I}}\right)I\,.
  \label{eq:eta_I}
\end{equation}
Thus, the corresponding prolongations simplify to
\begin{align}
\eta_{E}^{(1)}&=-\left(\dfrac{\chi_{k_{E}}}{k_{E}}\right)\dfrac{\mathrm{d}E}{\mathrm{d}t}\,.\label{eq:eta_E_prolong}\\
  \eta_{I}^{(1)}&=-\left(\dfrac{\chi_{k_{I}}}{k_{I}}\right)\dfrac{\mathrm{d}I}{\mathrm{d}t}\,.\label{eq:eta_I_prolong}
\end{align}
Given these simplifying assumptions and conditions, we can now assemble the linearised symmetry conditions.

Starting with the linearised symmetry condition for $E$, it is given by
\begin{equation}
  \begin{split}
    - \left(\beta \left(1 - \upsilon\right) I S - (\delta + \mu_{E}) E\right) \left(\dfrac{\chi_{k_{E}}}{k_{E}}\right) &= \beta \left(1 - \upsilon\right) I \eta_{S} - \beta I S \chi_{\upsilon} - \beta \left(1 - \upsilon\right) I S \left(\dfrac{\chi_{k_{I}}}{k_{I}}\right)\\
    &\quad+ \left(1 - \upsilon\right) I S \chi_{\beta} - E (\chi_{\delta} + \chi_{\mu_{E}}) + \left(\delta + \mu_{E}\right) E \left(\dfrac{\chi_{k_{E}}}{k_{E}}\right)\,.
  \end{split}
  \label{eq:lin_sym_E}
\end{equation}
Similarly, the linearised symmetry condition for $I$ is given by
\begin{equation}
  \begin{split}
    - \left(\beta \upsilon I S + \delta E - \mu_{I} I\right) \left(\dfrac{\chi_{k_{I}}}{k_{I}}\right) &= \beta \upsilon I \eta_{S} + \beta I S \chi_{\upsilon} - \delta E \left(\dfrac{\chi_{k_{E}}}{k_{E}}\right) + \upsilon I S \chi_{\beta} + E \chi_{\delta}\\
    &\quad- I \chi_{\mu_{I}} - \left(\beta \upsilon S - \mu_{I}\right) I \left(\dfrac{\chi_{k_{I}}}{k_{I}}\right)\,.
\end{split}
\label{eq:lin_sym_I}
\end{equation}
Next, we solve the linearised symmetry conditions for $E$ and $I$ in Eqs. \eqref{eq:lin_sym_E} and \eqref{eq:lin_sym_I}, respectively, for $\eta_S$ and equate them. After some simplifications, we obtain the following equation

\begin{equation}
\begin{split}
  &\quad- \beta k_{E} k_{I} I S \chi_{\upsilon} + \beta k_{E} \upsilon^{2} I S \chi_{k_{I}} - \beta k_{E} \upsilon I S \chi_{k_{I}} - \beta k_{I} \upsilon^{2} I S \chi_{k_{E}} + \beta k_{I} \upsilon I S \chi_{k_{E}} \\
  &\quad+ \delta k_{E} \upsilon E \chi_{k_{I}} - \delta k_{E} E \chi_{k_{I}} - \delta k_{I} \upsilon E \chi_{k_{E}} + \delta k_{I} E \chi_{k_{E}} - k_{E} k_{I} \upsilon E \chi_{\mu_{E}}\\
  &\quad- k_{E} k_{I} \upsilon I \chi_{\mu_{I}} - k_{E} k_{I} E \chi_{\delta} + k_{E} k_{I} I \chi_{\mu_{I}} = 0\,.
  \end{split}
  \label{eq:eta_S_equality}
\end{equation}
This equation depend on six unknown infinitesimals which we collect in the vector $\boldsymbol{\chi}_{0}\in\mathbb{R}^{6}$ given by
\begin{equation}
  \boldsymbol{\chi}_{0}=\begin{pmatrix}\chi_{\upsilon}\\\chi_{k_{I}}\\\chi_{k_{E}}\\\chi_{\mu_{E}}\\\chi_{\mu_{I}}\\\chi_{\delta}\end{pmatrix}\,.
  \label{eq:chi_vec_sub}
\end{equation}
Next, we split up Eq. \eqref{eq:eta_S_equality} with respect to iterated monomials of $\{S,E,I\}$. This procedure yields the following linear equations.

\begin{align}
I:&\quad- k_{E} k_{I} \upsilon \chi_{\mu_{I}} + k_{E} k_{I} \chi_{\mu_{I}}&=0\,,\label{eq:det_eq_0}\\
E:&\quad\delta k_{E} \upsilon \chi_{k_{I}} - \delta k_{E} \chi_{k_{I}} - \delta k_{I} \upsilon \chi_{k_{E}} + \delta k_{I} \chi_{k_{E}} - k_{E} k_{I} \upsilon \chi_{\mu_{E}} - k_{E} k_{I} \chi_{\delta}&=0\,,\label{eq:det_eq_1}\\
I S:&\quad- \beta k_{E} k_{I} \chi_{\upsilon} + \beta k_{E} \upsilon^{2} \chi_{k_{I}} - \beta k_{E} \upsilon \chi_{k_{I}} - \beta k_{I} \upsilon^{2} \chi_{k_{E}} + \beta kI \upsilon \chi_{k_{E}}&=0\,.\label{eq:det_eq_2}
\end{align}
The above linear system of equations can be expressed as the following matrix system $M\boldsymbol{\chi_{0}}=\boldsymbol{0}$
\begin{equation}
  \begin{split}
    &\underset{=M}{\underbrace{\begin{pmatrix}0 & 0 & 0 & 0 & - k_{E} k_{I} \upsilon + k_{E} k_{I} & 0\\0 & \delta k_{E} \upsilon - \delta k_{E} & - \delta k_{I} \upsilon + \delta k_{I} & - k_{E} k_{I} \upsilon & 0 & - k_{E} k_{I}\\- \beta k_{E} k_{I} & \beta k_{E} \upsilon^{2} - \beta k_{E} \upsilon & - \beta k_{I} \upsilon^{2} + \beta k_{I} \upsilon & 0 & 0 & 0\end{pmatrix}}}\underset{=\boldsymbol{\chi}_{0}}{\underbrace{\begin{pmatrix}\chi_{\upsilon}\\\chi_{k_{I}}\\\chi_{k_{E}}\\\chi_{\mu_{E}}\\\chi_{\mu_{I}}\\\chi_{\delta}\end{pmatrix}}}\\
    &=\underset{=\mathbf{0}}{\underbrace{\begin{pmatrix}0\\0\\0\end{pmatrix}}}\,,
   \end{split}
  \label{eq:mat_sys_1}
\end{equation}
where $M\in\mathbb{R}^{3\times{6}}$, $\boldsymbol{\chi}_{0}\in\mathbb{R}^{6}$ and $\mathbf{0}\in\mathbb{R}^{3}$. Clearly, the parameter infinitesimal vector $\boldsymbol{\chi}_{0}\in\mathbb{R}^{6}$ that solves this matrix equation is given by a linear combination of the basis vectors of the nullspace of the matrix $M$ which is denoted by $\mathcal{N}(M)$. This nullspace is three-dimensional and is given by
\begin{equation}
  \mathcal{N}(M)=\left\{ \begin{pmatrix}0\\\frac{k_{I}}{k_{E}}\\1\\0\\0\\0\end{pmatrix}, \  \begin{pmatrix}\frac{\upsilon^{2}}{\delta}\\\frac{k_{I} \upsilon}{\delta \upsilon - \delta}\\0\\1\\0\\0\end{pmatrix}, \  \begin{pmatrix}\frac{\upsilon}{\delta}\\\frac{k_{I}}{\delta \upsilon - \delta}\\0\\0\\0\\1\end{pmatrix}\right\}\,.
  \label{eq:nullspace_chi}
\end{equation}
Given three arbitrary coefficients $\alpha_{1},\alpha_{2},\alpha_{3}\in\mathbb{R}$, we have the following solution for six of the parameter infinitesimals:

\begin{equation}
  \begin{pmatrix}\chi_{\upsilon}\\\chi_{k_{I}}\\\chi_{k_{E}}\\\chi_{\mu_{E}}\\\chi_{\mu_{I}}\\\chi_{\delta}\end{pmatrix}=\begin{pmatrix}\upsilon \left(\alpha_{2} \upsilon + \alpha_{3}\right)\\\frac{k_{I} \left(\alpha_{1} \left(\upsilon - 1\right) + \alpha_{2} \upsilon + \alpha_{3}\right)}{\upsilon - 1}\\\alpha_{1} k_{E}\\\alpha_{2} \delta\\0\\\alpha_{3} \delta\end{pmatrix}\,.
\label{eq:chi_0_sol}
\end{equation}
From these calculations, we see immediately that $\chi_{\mu_{I}}=0$ and hence $\mu_{I}$ is locally structurally identifiable. After substituting these values into the equation for $\eta_{S}$, this infinitesimal is given by
\begin{equation}
  \eta_{S}=- \frac{a_{2} \beta \upsilon^{2} I S}{\beta \upsilon I - \beta I} + \frac{a_{2} \beta \upsilon I S}{\beta \upsilon I - \beta I} - \frac{a_{2} \delta E}{\beta \upsilon I - \beta I} - \frac{a_{3} \beta \upsilon I S}{\beta \upsilon I - \beta I} + \frac{a_{3} \beta I S}{\beta \upsilon I - \beta I} - \frac{a_{3} \delta E}{\beta \upsilon I - \beta I} - \frac{\upsilon I S \chi_{\beta}}{\beta \upsilon I - \beta I} + \frac{I S \chi_{\beta}}{\beta \upsilon I - \beta I}\,.
  \label{eq:eta_S}
\end{equation}
As a result of substituting the infinitesimal $\eta_{S}$ in Eq. \eqref{eq:eta_S} into its linearised symmetry condition, and then decomposing the resulting linearised symmetry condition with respect to iterated monomials of $\{S,E,I\}$, we obtain the following system of equations for the three remaining parameter infinitesimals $\chi_{\mu_{S}}$, $\chi_{\beta}$ and $\chi_{c}$ that are not part of the vector $\boldsymbol{\chi}_{0}$ in Eq. \eqref{eq:chi_vec_sub}:

\begin{align}
I^{2}:&\quad- \alpha_{2} \beta c \upsilon^{2} + \alpha_{2} \beta c \upsilon - \alpha_{3} \beta c \upsilon + \alpha_{3} \beta c - \beta \upsilon \chi_{c} + \beta \chi_{c} - c \upsilon \chi_{\beta} + c \chi_{\beta}&=0\,,\label{eq:det_eq_0_1}\\
E I:&\quad\alpha_{2} \delta^{2} + \alpha_{2} \delta \mu_{E} - \alpha_{2} \delta \mu_{I} - \alpha_{2} \delta \mu_{S} + \alpha_{3} \delta^{2} + \alpha_{3} \delta \mu_{E} - \alpha_{3} \delta \mu_{I} - \alpha_{3} \delta \mu_{S}&=0\,,\label{eq:det_eq_1_1}\\
E^{2}:&\quad\alpha_{2} \delta^{2} + \alpha_{3} \delta^{2}&=0\,,\label{eq:det_eq_2_1}\\
E I^{2}:&\quad- \alpha_{2} \beta \delta - \alpha_{3} \beta \delta&=0\,,\label{eq:det_eq_3_1}\\
I^{2} S:&\quad\alpha_{2} \beta \delta \upsilon - \alpha_{2} \beta \delta + \alpha_{3} \beta \delta \upsilon - \alpha_{3} \beta \delta + \beta \upsilon \chi_{\mu_{S}} - \beta \chi_{\mu_{S}}&=0\,,\label{eq:det_eq_4_1}\\
E I S:&\quad\alpha_{2} \beta \delta \upsilon + \alpha_{3} \beta \delta \upsilon&=0\,,\label{eq:det_eq_5_1}\\
I^{3} S:&\quad- \alpha_{1} \beta^{2} \upsilon + \alpha_{1} \beta^{2} - \alpha_{2} \beta^{2} \upsilon - \alpha_{3} \beta^{2} + \beta \upsilon \chi_{\beta} - \beta \chi_{\beta}&=0\,.\label{eq:det_eq_6_1}
\end{align}
We see from Eqs. \eqref{eq:det_eq_2_1}, \eqref{eq:det_eq_3_1} and \eqref{eq:det_eq_5_1} that
\begin{equation}
  \alpha_{2}=-\alpha_{3}\,.
  \label{eq:alpha_eq}
\end{equation}
By substituting this value, the above system is reduced to

\begin{align}
I^{2}:&\quad- \alpha_{2} \beta c \upsilon^{2} + 2 \alpha_{2} \beta c \upsilon - \alpha_{2} \beta c - \beta \upsilon \chi_{c} + \beta \chi_{c} - c \upsilon \chi_{\beta} + c \chi_{\beta}&=0\,,\label{eq:det_eq_new_0}\\
I^{2} S:&\quad\beta \left(\upsilon - 1\right) \chi_{\mu_{S}}&=0\,,\label{eq:det_eq_new_1}\\
I^{3} S:&\quad\beta \left(- \alpha_{1} \beta \upsilon + \alpha_{1} \beta - \alpha_{2} \beta \upsilon + \alpha_{2} \beta + \upsilon \chi_{\beta} - \chi_{\beta}\right)&=0\,.\label{eq:det_eq_new_2}
\end{align}
From Eq. \eqref{eq:det_eq_new_1}, we see that
\begin{equation}
  \chi_{\mu_{S}}=0\,,
  \label{eq:chi_muS}
\end{equation}
hence $\mu_{S}$ is locally structurally identifiable. Next, we solve Eq. \eqref{eq:det_eq_new_2} for $\chi_{\beta}$ which yields
\begin{equation}
  \chi_{\beta} = \beta \left(\alpha_{1} + \alpha_{2}\right)\,.
\label{eq:chi_beta}
\end{equation}
Lastly, we substitute the value for $\chi_{\beta}$ in Eq. \eqref{eq:chi_beta} into Eq. \eqref{eq:det_eq_new_0}, and solve the resulting equation for $\chi_{c}$ resulting in the following parameter infinitesimal
\begin{equation}
  \chi_{c} = - c \left(\alpha_{1} + \alpha_{2} \upsilon\right)\,.
  \label{eq:chi_c}
\end{equation}
All of these calculations result in the following vector for the parameter infinitesimals:
\begin{equation}
\begin{pmatrix}\chi_{c}\\\chi_{\beta}\\\chi_{\mu_{S}}\\\chi_{\mu_{E}}\\\chi_{\mu_{I}}\\\chi_{\delta}\\\chi_{\upsilon}\\\chi_{k_{E}}\\\chi_{k_{I}}\end{pmatrix}=\begin{pmatrix}- c \left(\alpha_{1} + \alpha_{2} \upsilon\right)\\\beta \left(\alpha_{1} + \alpha_{2}\right)\\0\\\alpha_{2} \delta\\0\\- \alpha_{2} \delta\\\alpha_{2} \upsilon \left(\upsilon - 1\right)\\\alpha_{1} k_{E}\\k_{I} \left(\alpha_{1} + \alpha_{2}\right)\end{pmatrix}\,,\label{eq:parameter_infinitesimals}
\end{equation}
which depend on the two arbitrary parameters $\alpha_{1}$ and $\alpha_{2}$. Moreover, this yields the following three infinitesimals for the states $S$, $E$ and $I$:
\begin{align}
  \eta_{S}&=-\left(\alpha_{1} + \alpha_{2} \upsilon\right) S,\label{eq:eta_S_final}\\
  \eta_{E}&=-\alpha_{1} E,\label{eq:eta_E_final}\\
\eta_{I}&=-(\alpha_{1}+\alpha_{2}) I.\label{eq:eta_I_final}  
\end{align}
In summary, this gives us the following family of infinitesimal generators of the Lie group:

\begin{equation}
  \begin{split}
    X &= -\left(\alpha_{1} + \alpha_{2} \upsilon\right) S\partial_{S} -\alpha_{1} E\partial_{E}-(\alpha_{1}+\alpha_{2}) I\partial_{I}\\
    &- c \left(\alpha_{1} + \alpha_{2} \upsilon\right)\partial_{c}+\beta \left(\alpha_{1} + \alpha_{2}\right)\partial_{\beta}+\alpha_{2} \delta\partial_{\mu_{E}}- \alpha_{2} \delta\partial_{\delta}\\
    &+\alpha_{2} \upsilon \left(\upsilon - 1\right)\partial_{\upsilon}+\alpha_{1}k_{E}\partial_{k_{E}}+(\alpha_{1}+\alpha_{2})k_{I}\partial_{k_{I}}\,,
\end{split}
  \label{eq:Lie_final_app}
\end{equation}
in terms of the two arbitrary parameters $\alpha_{1},\alpha_{2}\in\mathbb{R}$.

\subsection{Calculation of universal differential invariants of the SEI model}
The universal invariants above are independent of the arbitrary functions $\alpha_{1}$ and $\alpha_{2}$ respectively. Using the method of characteristics, we can integrate ODEs which are independent of these two functions and whose right hand sides are defined by the infinitesimals. The resulting integration constants correspond to the universal, since they are independent of the arbitrary functions $\alpha_{1}$ and $\alpha_{2}$, invariants which are the locally structurally identifiable and observable quantities.  

First, we see that none of the three states are locally structurally observable. But we can readily find observable quantitites depending on the states by calculating universal differential invariants. For the state $S$, we see that
\begin{equation}
  \dfrac{\mathrm{d}S}{\mathrm{d}c}=\dfrac{S}{c}\,,
  \label{eq:inv_S_eq}
\end{equation}
which is readily integrated to the following differential invariant
\begin{equation}
  I_{1}=\dfrac{S}{c}\,.
  \label{eq:inv_S}
\end{equation}
For the state $E$, we see that 
\begin{equation}
  \dfrac{\mathrm{d}E}{\mathrm{d}k_{E}}=-\dfrac{E}{k_{E}}\,,
  \label{eq:inv_E_eq}
\end{equation}
which is readily integrated to the following differential invariant
\begin{equation}
  I_{2}=k_{E}E\,.
  \label{eq:inv_E}
\end{equation}
Similarly, for the state $I$, we see that 
\begin{equation}
  \dfrac{\mathrm{d}I}{\mathrm{d}k_{I}}=-\dfrac{I}{k_{I}}\,,
  \label{eq:inv_I_eq}
\end{equation}
which is readily integrated to the following differential invariant
\begin{equation}
  I_{3}=k_{I}I\,.
  \label{eq:inv_I}
\end{equation}
Using the same approach, we next calculate the universal parameter invariants.

Since $\chi_{\mu_{S}}=0$ and $\chi_{\mu_{I}}=0$, the parameters $I_{4}=\mu_{S}$ and $I_{5}=\mu_{I}$ are locally structurally identifiable. Moreover, from the infinitesimals $\chi_{\mu_{E}}$ and $\chi_{\delta}$, we see that
\begin{equation}
  \dfrac{\mathrm{d}\mu_{E}}{\mathrm{d}\delta}=-1\,,
  \label{eq:inv_muE_delta_eq}
\end{equation}
which is integrated to give the following universal parameter invariant
\begin{equation}
  I_{6}=\mu_{E}+\delta\,.
  \label{eq:inv_muE_delta}
\end{equation}
From the infinitesimals $\chi_{k_{I}}$ and $\chi_{\beta}$, we see that  
\begin{equation}
  \dfrac{\mathrm{d}k_{I}}{\mathrm{d}\beta}=\dfrac{k_{I}}{\beta}\,,
  \label{eq:inv_kI_beta_eq}
\end{equation}
which is integrated to give the following universal parameter invariant
\begin{equation}
  I_{7}=\dfrac{k_{I}}{\beta}\,.
  \label{eq:inv_kI_beta}
\end{equation}
From the infinitesimals for $\chi_{\upsilon}$ and $\chi_{\delta}$, we see that  
\begin{equation}
  \dfrac{\mathrm{d}\delta}{\mathrm{d}\upsilon}=\dfrac{\delta}{\upsilon(1-\upsilon)}\,,
  \label{eq:inv_delta_upsilon_eq}
\end{equation}
which is integrated to give the following universal parameter invariant
\begin{equation}
  I_{8}=\delta\left(\dfrac{1-\upsilon}{\upsilon}\right)\,.
  \label{eq:inv_delta_upsilon}
\end{equation}
Next, we consider the combined infinitesimal for the product $c\upsilon$ which is given by
\begin{equation}
  \chi_{c\upsilon}=c\chi_{\upsilon}+\upsilon\chi_{c}=-(\alpha_{1}+\alpha_{2})c\upsilon\,,
  \label{eq:chi_product}
\end{equation}
and this product infinitesimal can be combined with $\chi_{\beta}$ which results in the equation
\begin{equation}
  \dfrac{\mathrm{d}(c\upsilon)}{\mathrm{d}\beta}=-\dfrac{(c\upsilon)}{\beta}\,.
  \label{eq:inv_c_upsilon_beta_eq}
\end{equation}
Integrating this equation yields the following universal parameter invariant
\begin{equation}
  I_{9}=\beta{c}\upsilon\,.
  \label{eq:inv_c_upsilon_beta}
\end{equation}
Finally, we consider the combined infinitesimal for the quotient $k_{E}/\delta$ which is given by
\begin{equation}
  \chi_{k_{E}/\delta}=\dfrac{\delta\chi_{k_{E}}-k_{E}\chi_{\delta}}{\delta^{2}}=(\alpha_{1}+\alpha_{2})\left(\dfrac{k_{E}}{\delta}\right)\,.
  \label{eq:kE_quotient}
\end{equation}
From this characteristic equation for the fraction $k_{E}/\delta$ and $\chi_{\beta}$, we see that
\begin{equation}
  \dfrac{\mathrm{d}\beta}{\mathrm{d}(k_{E}/\delta)}=\dfrac{\beta}{(k_{E}/\delta)}\,,
  \label{eq:inv_kE_delta_beta_eq}
\end{equation}
which is integrated to give the last universal parameter invariant;
\begin{equation}
  I_{10}=\dfrac{\beta\delta}{k_{E}}\,.
  \label{eq:inv_kE_delta_beta}
\end{equation}
 In summary, the ten universal invariants of the SEI model are given by
\begin{equation}
\begin{split}
    I_{1}&=\dfrac{S}{c}\,,\quad{I}_{2}=k_{E}E\,,\quad{I}_{3}=k_{I}I\,,\quad{I}_{4}=\mu_{S}\,,\quad{I}_{5}=\mu_{I}\,,\quad{I}_{6}=\mu_{E}+\delta\,,\\
    I_{7}&=\dfrac{k_{I}}{\beta}\,,\quad{I}_{8}=\delta\left(\dfrac{1-\upsilon}{\upsilon}\right)\,,\quad{I}_{9}=\beta{c}\upsilon\,,\quad{I}_{10}=\dfrac{\beta\delta}{k_{E}}\,.
\end{split}
    \label{eq:uni_inv_SEI_app}
\end{equation}


\end{document}